\documentclass[12pt]{article}
\usepackage{latexsym, amssymb}

\DeclareMathAlphabet\gothic{U}{euf}{m}{n}

\makeatletter
\def\eqnarray{\stepcounter{equation}\let\@currentlabel=\theequation
\global\@eqnswtrue
\tabskip\@centering\let\\=\@eqncr
$$\halign to \displaywidth\bgroup\hfil\global\@eqcnt\z@
  $\displaystyle\tabskip\z@{##}$&\global\@eqcnt\@ne
  \hfil$\displaystyle{{}##{}}$\hfil
  &\global\@eqcnt\tw@ $\displaystyle{##}$\hfil
  \tabskip\@centering&\llap{##}\tabskip\z@\cr}

\def\endeqnarray{\@@eqncr\egroup
      \global\advance\c@equation\m@ne$$\global\@ignoretrue}
\makeatother

\begin{document}
 \bibliographystyle{tom}

\newtheorem{lemma}{Lemma}[section]
\newtheorem{thm}[lemma]{Theorem}
\newtheorem{cor}[lemma]{Corollary}
\newtheorem{voorb}[lemma]{Example}
\newtheorem{rem}[lemma]{Remark}
\newtheorem{prop}[lemma]{Proposition}
\newtheorem{stat}[lemma]{{\hspace{-5pt}}}

\newenvironment{remarkn}{\begin{rem} \rm}{\end{rem}}
\newenvironment{exam}{\begin{voorb} \rm}{\end{voorb}}

\newcommand{\gota}{\gothic{a}}
\newcommand{\gotb}{\gothic{b}}
\newcommand{\gotc}{\gothic{c}}
\newcommand{\gote}{\gothic{e}}
\newcommand{\gotf}{\gothic{f}}
\newcommand{\gotg}{\gothic{g}}
\newcommand{\gothh}{\gothic{h}}
\newcommand{\gotk}{\gothic{k}}
\newcommand{\gotm}{\gothic{m}}
\newcommand{\gotn}{\gothic{n}}
\newcommand{\gotp}{\gothic{p}}
\newcommand{\gotq}{\gothic{q}}
\newcommand{\gotr}{\gothic{r}}
\newcommand{\gots}{\gothic{s}}
\newcommand{\gotu}{\gothic{u}}
\newcommand{\gotv}{\gothic{v}}
\newcommand{\gotw}{\gothic{w}}
\newcommand{\gotz}{\gothic{z}}
\newcommand{\gotG}{\gothic{G}}
\newcommand{\gotL}{\gothic{L}}
\newcommand{\gotS}{\gothic{S}}
\newcommand{\gotT}{\gothic{T}}

\newcounter{teller}
\renewcommand{\theteller}{\Roman{teller}}
\newenvironment{tabel}{\begin{list}%
{\rm \bf \Roman{teller}.\hfill}{\usecounter{teller} \leftmargin=1.1cm
\labelwidth=1.1cm \labelsep=0cm \parsep=0cm}
                      }{\end{list}}

\newcounter{tellerr}
\renewcommand{\thetellerr}{\roman{tellerr}}
\newenvironment{subtabel}{\begin{list}%
{\rm  \roman{tellerr}.\hfill}{\usecounter{tellerr} \leftmargin=1.1cm
\labelwidth=1.1cm \labelsep=0cm \parsep=0cm}
                         }{\end{list}}

\newcounter{proofstep}
\newcommand{\nextstep}{\refstepcounter{proofstep}\ruimte \par 
          \noindent{\bf Step \theproofstep} \hspace{5pt}}
\newcommand{\firststep}{\setcounter{proofstep}{0}\nextstep}

\newcommand{\Ni}{{\bf N}}
\newcommand{\Ri}{{\bf R}}
\newcommand{\Ci}{{\bf C}}
\newcommand{\Ti}{{\bf T}}
\newcommand{\Zi}{{\bf Z}}
\newcommand{\Fi}{{\bf F}}

\newcommand{\proof}{\mbox{\bf Proof} \hspace{5pt}} 
\newcommand{\remark}{\mbox{\bf Remark} \hspace{5pt}}
\newcommand{\ruimte}{\vskip10.0pt plus 4.0pt minus 6.0pt}

\newcommand{\ad}{{\mathop{\rm ad}}}
\newcommand{\Ad}{{\mathop{\rm Ad}}}
\newcommand{\Aut}{\mathop{\rm Aut}}
\newcommand{\arccot}{\mathop{\rm arccot}}
\newcommand{\diam}{\mathop{\rm diam}}
\newcommand{\divv}{\mathop{\rm div}}
\newcommand{\codim}{\mathop{\rm codim}}
\newcommand{\RRe}{\mathop{\rm Re}}
\newcommand{\IIm}{\mathop{\rm Im}}
\newcommand{\Tr}{{\mathop{\rm Tr}}}
\newcommand{\supp}{\mathop{\rm supp}}
\newcommand{\sgn}{\mathop{\rm sgn}}
\newcommand{\esssup}{\mathop{\rm ess\,sup}}
\newcommand{\Int}{\mathop{\rm Int}}
\newcommand{\Leibniz}{\mathop{\rm Leibniz}}
\newcommand{\lcm}{\mathop{\rm lcm}}
\newcommand{\mod}{\mathop{\rm mod}}
\newcommand{\spann}{\mathop{\rm span}}
\newcommand{\ubar}{\underline{\;}}
\newcommand{\one}{1\hspace{-4.5pt}1}

\hyphenation{groups}
\hyphenation{unitary}

\newcommand{\dm}{{d_{\mbox{\gothicss m}}}}
\newcommand{\dms}{{d_{\mbox{\gothics m}}}}
\newcommand{\hc}{\overline{H}}
\newcommand{\sodp}{ \{ 1,\ldots,d' \} }
\newcommand{\ccig}{C_c^\infty(G)}
\newcommand{\Cprime}[1]{C^{#1}{\hspace{0.8pt}}^\prime}
\newcommand{\CCprime}[1]{{\cal C}^{#1}{\hspace{0.8pt}}^\prime}
\newcommand{\Cdeltaprime}[1]{C_\Delta^{#1}{\hspace{0.5pt}}^\prime}

\newcommand{\Gdual}{\widetilde{G}}
\newcommand{\htozp}{H'_{2;1}\raisebox{10pt}[0pt][0pt]{\makebox[0pt]{\hspace{-28pt}$\scriptstyle\circ$}}}
\newcommand{\htozpu}{H^{\prime(u)}_{2;1}\raisebox{10pt}[0pt][0pt]{\makebox[0pt]{\hspace{-36pt}$\scriptstyle\circ$}}}

\newcommand{\htozpud}{H^{(u)}_{2;1}\raisebox{10pt}[0pt][0pt]{\makebox[0pt]{\hspace{-32pt}$\scriptstyle\circ$}}}
\newcommand{\htozpukd}{H^{(u_k)}_{2;1}\raisebox{10pt}[0pt][0pt]{\makebox[0pt]{\hspace{-41pt}$\scriptstyle\circ$}}}
\newcommand{\htozpund}{H^{(u_n)}_{2;1}\raisebox{10pt}[0pt][0pt]{\makebox[0pt]{\hspace{-41pt}$\scriptstyle\circ$}}}

\newcommand{\htozpH}{H^{\prime(H)}_{2;1}\raisebox{10pt}[0pt][0pt]{\makebox[0pt]{\hspace{-36pt}$\scriptstyle\circ$}}}
\newcommand{\htozpun}{H^{\prime(u_n)}_{2;1}\raisebox{10pt}[0pt][0pt]{\makebox[0pt]{\hspace{-45pt}$\scriptstyle\circ$}}}
\newcommand{\htoz}{H_{2;1}\raisebox{10pt}[0pt][0pt]{\makebox[0pt]{\hspace{-28pt}$\scriptstyle\circ$}}}
\newcommand{\twocases}{Suppose
   \begin{list}{}{\leftmargin=40mm \labelwidth=20mm \labelsep=5mm}
    \item[{\hfill\rm \bf either}\hfill] $H$ is  strongly elliptic 
    \item[{\hfill\rm \bf or}\hfill] $G$ is stratified and 
          $a_1,\ldots,a_{d'}$ is a basis for $\gotg_1$ 
     \newline
   in the stratification $(\gotg_m)_{m \in \{ 1,\ldots,r \} }$ of $\gotg$.
        \end{list}    }

\newcommand{\lstar}[1]{\mathbin{{}_{#1}\mskip-1mu*}}
\newcommand{\ldoublestar}[2]{\mathbin{{}_{#2}\mskip-1mu({}_{#1}\mskip-1mu*)}}
\newcommand{\sdp}[1]{\mathbin{{}_{#1}\mskip-2mu\times}}
\newcommand{\tfrac}[2]{{\textstyle \frac{#1}{#2} }}
\newcommand{\sptilde}{\hspace{0pt} \widetilde{\rule{0pt}{8pt} \hspace{7pt} }}

\newcommand{\cc}{{\cal C}}
\newcommand{\cd}{{\cal D}}
\newcommand{\cf}{{\cal F}}
\newcommand{\ch}{{\cal H}}
\newcommand{\ck}{{\cal K}}
\newcommand{\cl}{{\cal L}}
\newcommand{\cm}{{\cal M}}
\newcommand{\co}{{\cal O}}
\newcommand{\cs}{{\cal S}}
\newcommand{\ct}{{\cal T}}
\newcommand{\cx}{{\cal X}}
\newcommand{\cy}{{\cal Y}}
\newcommand{\cz}{{\cal Z}}

\newfont{\fontcmrten}{cmr10}
\newcommand{\slbrl}{\mbox{\fontcmrten (}}
\newcommand{\slbrr}{\mbox{\fontcmrten )}}

\newlength{\hightcharacter}
\newlength{\widthcharacter}
\newcommand{\covsup}[1]{\settowidth{\widthcharacter}{$#1$}\addtolength{\widthcharacter}{-0.15em}\settoheight{\hightcharacter}{$#1$}\addtolength{\hightcharacter}{0.1ex}#1\raisebox{\hightcharacter}[0pt][0pt]{\makebox[0pt]{\hspace{-\widthcharacter}$\scriptstyle\circ$}}}
\newcommand{\cov}[1]{\settowidth{\widthcharacter}{$#1$}\addtolength{\widthcharacter}{-0.15em}\settoheight{\hightcharacter}{$#1$}\addtolength{\hightcharacter}{0.1ex}#1\raisebox{\hightcharacter}{\makebox[0pt]{\hspace{-\widthcharacter}$\scriptstyle\circ$}}}
\newcommand{\scov}[1]{\settowidth{\widthcharacter}{$#1$}\addtolength{\widthcharacter}{-0.15em}\settoheight{\hightcharacter}{$#1$}\addtolength{\hightcharacter}{0.1ex}#1\raisebox{0.7\hightcharacter}{\makebox[0pt]{\hspace{-\widthcharacter}$\scriptstyle\circ$}}}
\newcommand{\gznp}{\settowidth{\widthcharacter}{$G$}\addtolength{\widthcharacter}{-0.15em}\settoheight{\hightcharacter}{$G$}\addtolength{\hightcharacter}{0.1ex}G\raisebox{\hightcharacter}[0pt][0pt]{\makebox[0pt]{\hspace{-\widthcharacter}$\scriptstyle\circ$}}_N'}
\newcommand{\sgznp}{{\settowidth{\widthcharacter}{$G$}\addtolength{\widthcharacter}{-0.15em}\settoheight{\hightcharacter}{$G$}\addtolength{\hightcharacter}{0.1ex}\raisebox{-5pt}{$\scriptstyle G\raisebox{0.7\hightcharacter}[0pt][0pt]{\makebox[0pt]{\hspace{-\widthcharacter}$\scriptstyle\circ$}}_N'$}}}

\thispagestyle{empty}

\begin{center}
{\Large{\bf Second-order operators with}}\\[2mm] 
{\Large{\bf degenerate coefficients}}  \\[2mm]
\large A.F.M. ter Elst$^1$, Derek W. Robinson$^2$, Adam Sikora$^3$
and Yueping Zhu$^4$
\end{center}

\vspace{5mm}

\begin{center}
{\bf Abstract}
\end{center}

\begin{list}{}{\leftmargin=1.8cm \rightmargin=1.8cm \listparindent=10mm 
   \parsep=0pt}
\item
We consider properties of second-order operators 
$H = -\sum^d_{i,j=1} \partial_i \, c_{ij} \, \partial_j$ on $\Ri^d$ with bounded
real symmetric measurable coefficients.
We assume that $C = (c_{ij}) \geq 0$ almost everywhere, but allow for the 
possibility that $C$ is singular.
We associate with $H$ a canonical self-adjoint viscosity operator $H_0$ and 
examine properties of the viscosity semigroup $S^{(0)}$ generated by $H_0$.
The semigroup extends to a positive contraction semigroup on the $L_p$-spaces
with $p \in [1,\infty]$.
We establish that 
it conserves probability, satisfies $L_2$~off-diagonal bounds
and  that the 
wave equation associated with $H_0$ has finite speed of propagation.
Nevertheless $S^{(0)}$ is not always strictly positive because separation
of the system can occur even for subelliptic operators. 
This demonstrates that subelliptic semigroups are   not ergodic in general
and  their  kernels are  neither  strictly positive nor H\"older continuous.
In particular one can construct examples for which both upper and lower Gaussian 
bounds fail even with coefficients in $C^{2-\varepsilon}(\Ri^d)$ with 
$\varepsilon > 0$.

\end{list}

\vspace{1cm}
\noindent
September 2005

\vspace{5mm}

\vspace{5mm}

\noindent
{\bf Home institutions:}    \\[3mm]
\begin{tabular}{@{}cl@{\hspace{10mm}}cl}
1. & Department of Mathematics  & 
  2. & Centre for Mathematics   \\
& \hspace{15mm} and Computing Science & 
  & \hspace{15mm} and its Applications  \\
& Eindhoven University of Technology & 
  & Mathematical Sciences Institute  \\
& P.O. Box 513 & 
  & Australian National University  \\
& 5600 MB Eindhoven & 
  & Canberra, ACT 0200  \\
& The Netherlands & 
  & Australia  \\[8mm]
3. & Department of Mathematical Sciences & 
 4. &  Department of Mathematics   \\
& New Mexico State University & 
  & Nantong University  \\
& P.O. Box 30001 & 
  & Nantong, 226007   \\
& Las Cruces & 
  & Jiangsu Province  \\
& NM 88003-8001, USA & 
  & P.R. China  \\
& {} & 
  & {} 
\end{tabular}

\newpage
\setcounter{page}{1}

\section{Introduction}\label{Sadt1}

Our intention is to investigate global properties of  second-order  operators $H$ 
with real measurable coefficients  on $\Ri^d$.
We consider  operators in divergence form formally given by
\begin{equation}
H=-\sum^d_{i,j=1}\partial_i\,c_{ij}\,\partial_j
\label{eadt1.1}
\end{equation}
where $\partial_i=\partial/\partial x_i$.
The coefficients $c_{ij}$ are assumed to be real $L_\infty$-functions
and the corresponding matrix $C=(c_{ij})$ is assumed to be symmetric and  positive-definite 
almost-everywhere.
Since  the classical work of Nash \cite{Nash} and De Giorgi \cite{DG}
the theory of such  operators is well developed under the additional hypothesis of 
strong ellipticity, i.e., the  assumption 
\begin{equation}
C\geq \mu \,I>0
\label{eadty1.2}
\end{equation}
almost-everywhere.
The principal result of this  theory is the local H{\"o}lder  continuity of weak solutions
of the associated elliptic and parabolic equations.
In Nash's approach  the H{\"o}lder
continuity of the elliptic solution is derived as a 
corollary of continuity of   the parabolic solution and the latter is established 
by an iterative argument from good upper and lower bounds
on the fundamental solution.
Aronson \cite{Aro} subsequently improved Nash's bounds and proved that the 
fundamental solution of the parabolic equation, the heat kernel, satisfies 
Gaussian upper and  lower bounds.
Specifically the kernel $K$  of the semigroup $S$ generated by $H$
is a symmetric  function over $\Ri^d\times\Ri^d$ satisfying bounds
\begin{equation}
a'\,G_{b';t}(x-y)\leq K_t(x\,;y)\leq a\,G_{b;t}(x-y)
\label{eo1.3}
\end{equation}
uniformly for  $x,y\in\Ri^d$ and $t>0$ where $G_{b;t}$ is the usual Gaussian function, 
$G_{b;t}(x)=t^{-d/2}e^{-b|x|^2 t^{-1}}$,
and $a,a',b,b'>0$.
(Background material on the Nash--De Giorgi theory can be found in the books and reviews
 \cite{Dav2} \cite{Gia1} \cite{Gia} \cite{Stroock1} \cite{Stroock2}.
The derivation of H\"older continuity from the Aronson bounds is well described in \cite{FaS}
and a clear  statement of the equivalence of estimates for elliptic and parabolic solutions
is given in \cite{Aus1}.)

In contrast to the Nash--De Giorgi theory we examine  operators for which the strong ellipticity
assumption (\ref{eadty1.2}) is not satisfied.
Part of our work requires nothing other than the ellipticity property $C\geq0$ but we also
analyze  operators which satisfy a  condition of subellipticity.
The strong ellipticity bound (\ref{eadty1.2}) on the coefficients  is equivalent to the operator bound
\begin{equation}
H\geq\mu \,\Delta
\label{eadty1.3}
\end{equation}
on $L_2(\Ri^d)$ where $\Delta=-\sum^d_{i=1}\partial_i^2$ is the usual self-adjoint Laplacian.
The subelliptic condition which we consider is 
\begin{equation}
H\geq \mu\,\Delta^\gamma-\nu\,I
\label{eadt1.4}
\end{equation}
where $\mu>0$, $\nu\geq 0$ and $\gamma\in \langle0,1]$.
This subellipticity condition  first arose in H\"ormander's  work \cite{Hor1} 
on the characterization of hypoelliptic operators as
sums of squares of $C^\infty$-vector fields satisfying a fixed rank condition. 
This work was extended by Rothschild and Stein \cite{RS} and the relation between the rank~$r$ 
of  the vector fields and the order $\gamma$ was clarified.
In fact $\gamma=1/r\in \langle0,1/2] \cup \{ 1 \} $.
Subsequently Fefferman and Phong \cite{FP} (see also \cite{FSC} \cite{SC} and \cite{OR}) 
analyzed operators with smooth coefficients satisfying the inequality without 
assuming  they could be expressed
as sums of squares of vector fields.
They established that the subellipticity condition could be characterized by properties
of the intrinsic geometry.
All this analysis was of a local nature.
Later Kusuoka and Stroock \cite{KuS} examined global properties of operators of the form
(\ref{eadt1.1}) under various assumptions on the local geometry and   positivity of the 
corresponding semigroup kernels (see, for example, Theorems~(2.6), (3.1) and (3.9) of \cite{KuS}).
These  results could then be applied to sums of squares of vector fields satisfying a uniform 
version of H\"ormander's rank condition (see \cite{KuS} Theorems~(3.20) and (3.24)).
It is notable that many of the estimates of Fefferman--Phong and Kusuoka--Stroock only depend
on the $C^2$-norm of the coefficients.
One can,  however, establish a broad range of examples for which the Fefferman--Phong 
characterization of subellipticity and the  lower bounds on the kernel
fail if the  coefficients are not in $C^2(\Ri^d)$.
In particular the global behaviour of the kernel is quite different to the smooth situation.

The results of Fefferman--Phong indicate that the local behaviour is governed by the  intrinsic
geometry associated with the subelliptic operator.
The Kusuoka--Stroock philosophy, explained in the introduction to \cite{KuS}, is based on the idea
that the detailed geometry is blurred with passing time  and that the semigroup kernel should resemble
the standard Gaussian $G_{b,t}$ for large time.
Our results establish that this is not the case for a large class of subelliptic operators
whose coefficients are less than twice differentiable.
Local properties often persist and dictate the global behaviour.

The phenomenon which distinguishes between general subelliptic operators with measurable coefficients
and those of the H\"ormander type is the possibility  of separation.
For example, in one-dimension the operator $H=-d\,c\,d$ satisfies the subelliptic condition (\ref{eadt1.4})
if~$c$ has an isolated zero $c(x) \asymp x^{2(1-\gamma)}$ as $x\to 0$ with 
$\gamma\in\langle 0,1/2 \rangle$.
Nevertheless
$H$ separates into a direct sum of two operators acting on the left and right half-lines, respectively.
Then the corresponding kernel cannot be strictly positive nor uniformly continuous.
More complicated separation phenomena occur if $c$ has several zeros or for operators in higher dimensions.

The theory of elliptic operators in divergence form, and in non-divergence form, has a long
and complex history.
A partial perspective on modern aspects of the fundamental theory can be obtained from the
books \cite{Fri} \cite{Gia} \cite{GT} \cite{Hor4} \cite{Hor5} \cite{Hor6} \cite{Hor7}
\cite{Tay3} \cite{Tre1} and references therein.
Probabilistic methods and stochastic analysis have  been applied to the analysis of 
elliptic operators and relevant information can be found in \cite{Stroock3} \cite{SV1}.
More recently the theory has been extended to the setting of Dirichlet spaces 
(see, for example, \cite{BM} \cite{Stu3} \cite{Stu4} \cite{Stu5}).

Although the theory of strongly elliptic operators 
in divergence form  is now well understood and systematically developed
the same cannot be said of the theory of degenerate elliptic operators.
Despite much interest in degenerate operators
(see, for example, \cite{BM} \cite{Fra} \cite{FKS} \cite{FL}  \cite{FLW} \cite{LaM} 
 \cite{MV1} \cite{Mus} \cite{Tru2} \cite{Varo10} and references therein)  there is no
commonly accepted definition of degeneracy.
Various conditions of positivity, integrability and regularity of the lowest eigenvalue
of the coefficient matrix $C$  have been proposed and studied as measures of degeneracy.
The main aim of many of the investigations have been to prove H\"older continuity 
of solutions or to derive Poincar\'e or Sobolev style inequalities,
properties analogous to those of strongly
elliptic operators.
Our results are of a different nature.
We examine situations in which the corresponding heat kernels are not even continuous
(see Section~\ref{Sadty5}).
Therefore many of the regularity conditions analyzed in the previous works are not
satisfied.
We stress, however, that one can nevertheless obtain many positive results for operators 
 with irregular and degenerate coefficients, e.g.,  $L_2$ off-diagonal bounds, 
finite speed of propagation of the corresponding wave equation
and lower bounds for the kernel of  high powers of the resolvent.

\section{Elliptic operators}\label{Sadty2}
The first problem in the analysis of the  elliptic operators  (\ref{eadt1.1})
is    the rigorous definition of $H$ as a 
positive self-adjoint operator on $L_2(\Ri^d)$.
This is a delicate problem for degenerate operators although the 
delicacies are often overlooked. 
The usual approach in operator theory is by quadratic forms.
First one introduces the elliptic form 
\begin{equation}
h(\varphi)=\sum^d_{i,j=1}(\partial_i\varphi,c_{ij}\partial_j\varphi)
\label{eadty2.1}
\end{equation}
where $D(h)=\bigcap^d_{i=1}D(\partial_i)=L_{2;1}(\Ri^d)=D(\Delta^{1/2})$.
Then $h$ is  positive, symmetric and densely-defined.
Therefore if $h$ is closed there is a uniquely defined positive  self-adjoint operator $H$ 
such that $D(H)\subseteq D(h)$,  $D(H^{1/2})=D(h)$ and  $h(\varphi)=\|H^{1/2}\varphi\|_2^2$ (see, for example,
\cite{Kat1} Section~VI.2).
Alternatively if $h$ is closable then one can define $H$ in a similar manner through the
closure of $h$.
Therefore the first onus of any careful investigation is to establish closure properties
 of the quadratic form $h$.

If the  coefficients of the operator satisfy the strong ellipticity assumption 
(\ref{eadty1.2})  it is easy to deduce that $h$ is closed. 
Then
\begin{equation}
\|C\|\,l(\varphi)\geq h(\varphi)\geq \mu\,l(\varphi)
\label{eadty2.2}
\end{equation}
for all $\varphi\in D(h)$ where $\|C\|$ is the essential supremum of the norms of the 
matrices $C(x)$ and $l$ is the  quadratic form of the Laplacian, 
$l(\varphi)=\|\Delta^{1/2}\varphi\|_2^2$.
It follows immediately that $h$ is closed.
More is true.

\begin{prop}\label{padty2.1}
The form $(\ref{eadty2.1})$ is closed on $D(h)=L_{2;1}(\Ri^d)$  if and only if 
$h$ is strongly elliptic, i.e., $h(\varphi)\geq \mu\,l(\varphi)$ for some $\mu>0$ and all
$\varphi\in D(h)$.
\end{prop}
\proof\ Strong ellipticity implies that $h$ is closed by the foregoing comparison 
(\ref{eadty2.2}).
Conversely, if the form is closed then $D(h)$ is a Hilbert
space under the norm $|||\varphi|||^2 = h(\varphi) + \|\varphi\|_2^2$, 
by \cite{Kat1}, Theorem VI.1.11.
Alternatively, $D(\Delta^{1/2})$ is a Hilbert space with the graph norm.
Moreover, $D(\Delta^{1/2}) = L_{2;1}(\Ri^d) = D(h)$ as sets and both the spaces
$D(h)$ and $D(\Delta^{1/2})$ are continuously embedded in $L_2(\Ri^d)$. 
Hence,  by the closed graph theorem, there is a $\mu > 0$ such that 
$|||\varphi|||^2 \geq \mu \|(I + \Delta)^{1/2} \varphi\|_2^2$ for all 
$\varphi \in L_{2;1}(\Ri^d)$. 
Then it follows that 
\[
h(\varphi)\geq \mu\,l(\varphi)-(1-\mu)\,\|\varphi\|_2^2
\]
for all $\varphi\in D(h)$.
Now one can evaluate this inequality with $\varphi$ replaced by $\varphi_k$ where
$\varphi_k(x)=e^{ikx.\xi}\varphi(x)$ with $k\in\Ri$, $\xi\in\Ri^d$ 
and $\varphi\in C_c^\infty(\Ri^d)$.
Then using $\|\varphi_k\|_2=\|\varphi\|_2$ one calculates that
\[
\int_{\Ri^d}dx\,|\varphi(x)|^2\,(\xi,C(x)\xi)=
\lim_{k\to\infty}k^{-2}h(\varphi_k)\geq 
\lim_{k\to\infty}k^{-2}\mu\,l(\varphi_k)=\mu\,|\xi|^2\,\|\varphi\|_2^2
\;\;\;.
\]
Therefore $C\geq\mu\,I$ almost everywhere and $h \geq \mu \, l$.\hfill$\Box$

\ruimte

\begin{remarkn}\label{radty2.1}
This argument establishes that the subellipticity condition 
(\ref{eadt1.4}) with $\gamma=1$
is equivalent to the strong ellipticity condition (\ref{eadty1.3}),
i.e., if $\gamma=1$ one can choose $\nu=0$.
\end{remarkn}

Other criteria for the form (\ref{eadty2.1}) to be closed are given in \cite{ERZ1}.
In general the form $h$ is not closable. 
Nevertheless there are various useful criteria for closability.

First the above comparison argument gives a general criterion for closability.
If $k$ is a closable form with $D(k)=D(h)$ and one has estimates
\[
a_1\,k(\varphi) \leq h(\varphi) \leq a_2\,k(\varphi) 
\]
for some $a_1,a_2>0$ and all $\varphi\in D(h)$ then $h$ is closable.
This reasoning  can be applied to some degenerate operators.
For example, if $h$ is the form of an elliptic operator and $\mu_m(x), \mu_M(x)$ 
denote the smallest and largest eigenvalues of the coefficient matrix $C(x)$ for all $x \in \Ri^d$
and if
$\mu_M\leq a\,\mu_m$ almost everywhere for some constant $a>0$,  
then the closability of $h$ is equivalent
to closability of the form $k(\varphi)=\sum^d_{i=1}(\partial_i\varphi,\mu_m\partial_i\varphi)$
with $D(k)=D(h)$.

Secondly, if the coefficients $c_{ij}$ are once-continuously differentiable then $h$ is closable
since  $H$ can be  identified as a symmetric operator 
$H=-\sum^d_{i,j=1}c_{ij}\partial_i\partial_j
-\sum^d_{j=1}(\sum^d_{i=1}\partial_ic_{ij})\partial_j$
with domain $D(H)= L_{2;2}(\Ri^d)=\bigcap^d_{i,j=1}D(\partial_i\partial_j)$.
Then $h(\varphi)=(\varphi, H\varphi)$ for all $\varphi\in D(H)$ 
and the form $h$ is closable by the Friederich's extension method 
(see, \cite{Kat1}, Section VI.2.3).

Thirdly, in one-dimension a complete characterization of  closability is given in \cite{FOT},
pages~105--107.
This observation allows one to construct examples of non-closable elliptic~$h$ 
in higher dimensions.
Moreover the argument gives a sufficient condition for closability in higher dimensions
which covers many situations of degeneracy.

\begin{prop}\label{padty2.2}
Let $\mu_m(x)$ be the smallest eigenvalue of the coefficient matrix
$C(x)$ for all $x \in \Ri^d$.
Suppose for almost every $x \in \Ri^d$ there exists a neighbourhood $U$ of $x$ 
such that $\mu_m > 0$ almost everywhere on $U$ and $\int_U \mu_m^{-1} < \infty$.
Then the form  $h$ is closable.
\end{prop}
\proof\
See \cite{MR}, Section~II.2b. \hfill$\Box$

\ruimte

Many of the examples which we subsequently consider are covered by the next corollary.

\begin{cor} \label{cadty240} 
If $\mu_m$ is continuous and has a discrete set of zeros then the form $h$ is closable.
\end{cor}

In the general situation, when it is unclear if $h$ is closable,
 we adopt a different approach to the definition of the elliptic operator.
We will define it by an approximation method akin to the viscosity method of partial 
differential equations.

 Define $h_\varepsilon$ for each  $\varepsilon\in\langle0,1]$ by 
$D(h_\varepsilon)=D(h)=D(l)$ and
\[
h_\varepsilon(\varphi)=h(\varphi)+\varepsilon\,l(\varphi) 
\;\;\;.
\]
Then $h_\varepsilon$ is the closed form associated with the strongly elliptic 
operator with coefficients $C_\varepsilon = C+\varepsilon I$.
But  $\varepsilon\mapsto h_\varepsilon(\varphi)$ decreases monotonically as 
$\varepsilon$ decreases
for each $\varphi\in D(h)$.
Therefore it  follows from a result of Kato, \cite{Kat1} Theorem~VIII.3.11, 
that the $H_\varepsilon$
converge in the strong resolvent sense, as $\varepsilon\to0$, 
to a positive self-adjoint operator
$H_0$ which we  refer to as the {\bf viscosity operator} with coefficients 
$C=(c_{ij})$.
This procedure gives a   precise meaning to the  operator $H$ formally given by 
(\ref{eadt1.1}).

In the following we frequently have to compare two forms and two self-adjoint operators.
If $k_1$ and $k_2$ are two  symmetric forms with domains $D(k_1)$ and $D(k_2)$ in the same 
Hilbert space, then we write $k_1 \leq k_2$ if $D(k_1) \supseteq D(k_2)$ and
$k_1(\varphi) \leq k_2(\varphi)$ for all $\varphi \in D(k_2)$.
Next, there is a one-one correspondence between lower bounded self-adjoint
operators and closed, densely defined, symmetric, lower bounded quadratic forms.
Hence, if, in addition, $k_1$ and $k_2$ are lower bounded, closed and densely defined
and $K_1$ and $K_2$ are the associated self-adjoint operators then we write 
$K_1 \leq K_2$ in the sense of quadratic forms if $k_1 \leq k_2$.

Let  $h_0$ denote the form associated with $H_0$,
i.e., $D(h_0)=D(H_0^{1/2})$ and $h_0(\varphi)=\|H_0^{1/2}\varphi\|_2^2$.
Although the construction of the form $h_0$ might appear arbitrary it does have
an interesting property of universality.

\begin{prop} \label{padty215}
The following are valid.
\mbox{}
\begin{tabel} 
\item \label{padty215-1}
The viscosity form $h_0$ is the largest positive, symmetric, closed,
quadratic form $k$ with $k \leq h$.
\item \label{padty215-2}
$h_0(\varphi)=  h(\varphi)$ for all $\varphi\in D(h)$ if and only if
$h$ is closable and then $h_0=\overline h$, the closure of $h$.
\end{tabel}
\end{prop}
\proof\
Simon \cite{bSim5} defines the regular part of a general positive symmetric densely-defined quadratic
form as the largest closable symmetric quadratic form $k$ with $k \leq h$.
Therefore, with this terminology, the first statement states that $h_0$ is the closure
of the regular part of $h$.
Then the Statement~\ref{padty215-1} follows directly from Theorem~3.2 of \cite{bSim5}.

If, however,  $h_0(\varphi)=h(\varphi)$ for all 
$\varphi\in D(h)$, then $h_0$ is a closed extension of $h$. 
Hence $h$ is closable.
Conversely if $h$ is closable then $h$ equals the regular part of $h$
and $h_0=\overline h$ by Statement~\ref{padty215-1}.\hfill$\Box$

\ruimte

One implication of Statement~\ref{padty215-1} of Proposition~\ref{padty215}
is that $h_0$ is independent of the particular
approximation technique we have used, i.e., the addition of a small multiple of the
Laplacian.
The same limit would be obtained if one were to use a multiple of the square of the 
Laplacian.
The latter would correspond more closely to a viscosity term.

Quadratic forms play a significant role in convex analysis \cite{ET} and convergence
theory \cite{Bra} \cite{DalM} but the emphasis is rather different to that of operator theory.
In these areas of application properties of lower semicontinuity are important.
If $h$ is a positive symmetric quadratic form defined on a dense subspace $D(h)$ of the Hilbert
space $\cal H$ and if one extends $h$ to $\cal H$ by setting $h(\varphi)=\infty$ if 
$\varphi\notin D(h)$ then $h$ is closed if and only if 
the extension is lower semicontinuous \cite{bSim4} \cite{Kat1},  Lemma VIII.3.14a.
In general the lower semicontinuous envelope of the extension, which is variously called the
lower semi continuous regularization of $h$ \cite{ET} page~10 or the  relaxed form \cite{DalM} page~28,
determines a closed quadratic form.
The latter is the closure of the form $h$, if $h$ is closable, or is the closure of 
the regular part, in
Simon's terminology, if $h$ is not closable.
The regularization, or relaxation, has been used in
a variety of applications to nonlinear phenomena and discontinuous media
(see, for example, \cite{Bra} \cite {ET} \cite{Jos} \cite{DalM}  \cite{Mosco}
and references therein).
Mosco gives examples, on pages 414--416,   of relaxed forms 
which can be traced back to the basic example given in Beurling and Deny's paper \cite{BeD}.

\smallskip

The viscosity operator $H_0$ generates a self-adjoint contraction 
 semigroup $S^{(0)}$ on $L_2(\Ri^d)$.
Then since $H_0$ is defined as the strong resolvent limit of the strongly elliptic
 operators $H_\varepsilon$ associated with the closed forms $h_\varepsilon$
it follows that the semigroup $S^{(0)}$ is the strong limit of
 the self-adjoint contraction semigroups  $S^{(\varepsilon)}$ 
generated by the $H_\varepsilon$.
In particular  $S^{(\varepsilon)}_t$ converges strongly, on $L_2(\Ri^d)$, to $S^{(0)}_t$ and the convergence
is uniform for $t$ in finite intervals.
Note that each  $h_\varepsilon$ is a Dirichlet form, i.e., it satisfies 
the Beurling--Deny criteria 
(see, for example, \cite{RS4}, Appendix to Section~XIII.12, or \cite{Dav2}, Section~1.3).
Specifically a positive, symmetric, closed, quadratic form $h$ on $L_2$ is a Dirichlet form
if it satisfies the  following two conditions:
\begin{enumerate}
\item
$\varphi\in D(h)$ implies $|\varphi|\in D(h)$ and
$h(|\varphi|)\leq h(\varphi)$,
\item
$\varphi\in D(h)$ implies $\varphi\wedge\one\in D(h)$ and
$h(\varphi\wedge\one)\leq h(\varphi)$.
\end{enumerate}
(For a full description of the theory of Dirichlet forms 
see \cite{BH} \cite{FOT} \cite{Sil}.)

The primary result of the Beurling--Deny theory is that $h$ is a Dirichlet form if
and only if the semigroup
$S$ generated by the corresponding operator 
$H$ on $L_2$ is positivity preserving  and extends from $L_2\cap L_p$
to a contraction semigroup on $L_p$ for all $p\in[1,\infty]$.
Therefore the  semigroups $S^{(\varepsilon)}$ are positivity preserving 
and extend to  positivity preserving contraction semigroups, also denoted by $S^{(\varepsilon)}$, 
on each of the spaces  $L_p(\Ri^d)$ with $p\in[1,\infty]$.
Since $S^{(\varepsilon)}$ converges strongly to $S^{(0)}$ on  $L_2(\Ri^d)$ it follows
that  $S^{(0)}$ is positivity preserving. 
It also  extends to a contraction 
semigroup on the $L_p$-spaces by observing that 
\[
|(\varphi,S^{(\varepsilon)}_t\psi)|\leq \|\varphi\|_p\,\|\psi\|_q
\]
for all $\varphi\in L_2\cap L_p$ and $\psi\in L_2\cap L_q$ where $p$ and $q$ are conjugate 
exponents.
Then   similar estimates  follow for $S^{(0)}_t$ by taking the limit $\varepsilon\to0$.
Therefore  $S^{(0)}$ extends to a contraction semigroup on all the $L_p$-spaces
by a density argument.
The resulting extensions are obviously positivity preserving and so  $h_0$ must be a  
Dirichlet form.

Finally we note that as a consequence of  positivity and contractivity  the 
viscosity semigroup $S^{(0)}$ satisfies the Markov property
\begin{equation}
0\leq S^{(0)}_t\one\leq \one
\label{eeo2.10}
\end{equation}
for all $t>0$ on $L_\infty(\Ri^d)$.
In the next section we will, however, prove that $S^{(0)}_t\one= \one$.
This stronger property is often referred to as conservation of probability or 
stochastic completeness.
It is the property that motivated the work of Gaffney \cite{Gaf}.

\section{$L_2$ off-diagonal estimates}\label{Sadty3}

One may associate with  the coefficients $C = (c_{ij})$ a
`distance' $d_C \colon \Ri^d \times \Ri^d \to [0,\infty]$  by setting
\begin{equation}
d_C(x\,;y)=\sup_{\psi\in{\cal D}}|\psi(x)-\psi(y)|
\label{eadty3.1}
\end{equation}
for all $x,y\in\Ri^d$, where
\begin{equation}
{\cal D}
=\{\psi\in C_c^\infty(\Ri^d) : \psi \mbox{ real and }
      \Big\|\sum^d_{i,j=1}c_{ij} \, (\partial_i\psi) \, (\partial_j\psi)\Big\|_\infty\leq 1 \}
\;\;\;.
\label{eadty3.2}
\end{equation}
If $C$ is strongly elliptic then 
it follows from the bounds $\|C\| \,I\geq C\geq \mu\,I$ that 
\[
\|C\|^{-1/2}\,|x-y|\leq d_C(x\,;y)\leq \mu^{-1/2}\,|x-y|
\]
for all $x,y\in \Ri^d$, i.e., $d_C$ is a proper distance and it is 
equivalent to the Euclidean distance.
For degenerate operators,  however, $d_C$ is a pseudodistance,
i.e., it has the metric properties of a distance but it  can take the value infinity.
Nevertheless, for brevity we will refer to it as a distance. 

There are a variety of other methods of associating a distance with $C$ especially  if the 
coefficients are continuous. 
Then one may adopt one of several equivalent `shortest path' definitions (see \cite{JSC1} for
a survey and comparison of various possibilities).
A definition of the foregoing nature was introduced by Biroli and Mosco \cite{BM} in the
general context of Dirichlet forms and this was crucial for the extension of 
many concepts of  elliptic operator theory to this setting \cite{Stu4} \cite{Stu5}.
In the case of degenerate $C$  it is not evident that (\ref{eadty3.1}) is the most appropriate 
definition (see Section~\ref{Sadty5}) but it is adequate for many purposes.

Our immediate purpose is to examine a general type of Gaussian bound on $L_2(\Ri^d)$
which  originated in the work of  Gaffney \cite{Gaf}.
Bounds of this type have subsequently been considered by various authors (see, for example, 
\cite{Aus2} \cite{CGT}  \cite{Dav12} \cite{Gri3}  \cite{Stu4} \cite{Stu2}).
The bounds are variously called integrated Gaussian estimates, an integrated
maximum principle or $L_2$ off-diagonal bounds.

For all $x \in \Ri^d$ and $r > 0$ set  $B_C(x\,;r)=\{y\in X: \,d_C(x\,;y)<r \}$.
In the sequel we fix $x_1,x_2 \in \Ri^d$ and $r_1,r_2 > 0$
and consistently use the notation $B_1=B_C(x_1\,;r_1)$ and 
$B_2=B_C(x_2\,;r_2)$ for balls and  set
\[
\tilde{d}_C(B_1\,;B_2)=(d_C(x_1\,;x_2)-r_1-r_2)\vee 0
\;\;\;.
\]
Note that it follows from the triangle inequality that
\[
\tilde{d}_C(B_1\,;B_2)\leq d_C(B_1\,;B_2)
= \inf_{x\in B_1}\inf_{y\in B_2}d_C(x\,;y)
\;\;\;.
\]
The subsequent  proof of  $L_2$ off-diagonal bounds for the  viscosity semigroup $S^{(0)}$
follows the  arguments of Davies \cite{Dav12}.
Care has to be taken since the distance $d_C$ 
can take the value infinity.
We adopt the convention $e^{-\infty} = 0$.

\begin{prop} \label{padty311}
The viscosity semigroup $S^{(0)}$ satisfies
\begin{equation}
|(\varphi_1,S^{(0)}_t\varphi_2)|\leq 
e^{-\tilde d_C(B_1;B_2)^2/(4t)}\, \|\varphi_1\|_2 \, \|\varphi_2\|_2
\label{etadt201;1}
\end{equation}
 $\varphi_1 \in L_2(B_1)$, $\varphi_2 \in L_2(B_2)$ and $t > 0$.
\end{prop}
\proof\
If $d_C(x_1\,;x_2)-r_1-r_2 \leq 0$ then (\ref{etadt201;1})  follows from the 
contractivity of $S^{(0)}$, so we may assume that 
$d_C(x_1\,;x_2)-r_1-r_2 > 0$.
In particular $d_C(x_1\,;x_2) \in \langle0,\infty]$.
Let $r \in \langle0,\infty\rangle$ and suppose that $r < d_C(x_1\,;x_2)$
and $r - r_1 - r_2 > 0$.
By definition of $d_C$ there exists a $\psi \in {\cal D}$ such that $\psi(x_2) - \psi(x_1) > r$.

Consider  the bounded multiplication operator $U_\rho$ 
defined by
$U_\rho \varphi = e^{- \rho\psi} \varphi$ for all $\rho \in \Ri$.
If $\varphi\in D(h)$  it follows  that $U_\rho\varphi\in D(h)$ and 
$\partial_iU_\rho\varphi=U_\rho(\partial_i-(\partial_i\psi))\varphi$.
Moreover, if $\varphi$ is  real
\begin{eqnarray*}
\sum^d_{i,j=1}(\partial_iU_\rho\varphi, c_{ij}\partial_jU_\rho^{-1}\varphi)
&=& h(\varphi)-\rho^2\sum^d_{i,j=1}
((\partial_i\psi)\varphi, c_{ij}(\partial_j\psi)\varphi)\\[5pt]
&\geq& h(\varphi)-\rho^2\,\|\varphi\|_2^2
\end{eqnarray*}
since the terms linear in $\rho$ cancel by reality and  symmetry.
Similarly,  if $H_\varepsilon$ are  the strongly elliptic approximants,
  with the coefficients $c_{ij}+\varepsilon\delta_{ij}$, to the
viscosity operator $H_0$  one has bounds
\[
h_\varepsilon(U_\rho\varphi,  U_\rho^{-1}\varphi)
\geq - \rho^2(1+\varepsilon \|\nabla\psi\|_\infty^2) \|\varphi\|_2^2
\]
for all real $\varphi\in D(h)=L_{2;1}$.
Therefore, using an obvious differential inequality, one deduces that
$\|U_\rho S^{(\varepsilon)}_t U_\rho^{-1}\|_{2\to2}\leq
e^{\rho^2(1+\varepsilon \|\nabla\psi\|_\infty^2) t}$. 
Then  by taking the the
limit $\varepsilon \downarrow 0$ one concludes that
\[
\|U_\rho S^{(0)}_t U_\rho^{-1}\|_{2\to2}\leq e^{\rho^2 t}
\]
for all $\rho \in \Ri$ and $t>0$.
The estimate is initially valid on the real $L_2$-space and then by
polarization
on the complex space.

Next, if $x \in B_1$ then $\psi(x) - \psi(x_1) \leq d_C(x\,;x_1) < r_1$.
So $\rho \psi(x) \leq \rho(\psi(x_1) + r_1)$ for all $\rho > 0$ and 
\[
\|U_\rho^{-1} \varphi_1\|_2
= \|e^{\rho \psi} \varphi_1\|_2
\leq e^{\rho(\psi(x_1)+ r_1)} \, \|\varphi_1\|_2
\;\;\;.
\]
Alternatively, if $x \in B_2$ then $|\psi(x) - \psi(x_2)| \leq d_C(x\,;x_2) < r_2$.
So 
\[
\psi(x) - \psi(x_1)
= \psi(x_2) - \psi(x_1) - \left( \psi(x_2) - \psi(x) \right)
> r - r_2
\]
and $\rho \psi(x) \geq \rho(\psi(x_1) + r - r_2)$ for all $\rho > 0$.
Therefore
\[
\|U_\rho \varphi_2\|_2
= \|e^{- \rho \psi} \varphi_2\|_2
\leq e^{-\rho(\psi(x_1)+r-r_2)} \, \|\varphi_2\|_2  
\;\;\; .  \]
Combining these estimates one deduces that 
\begin{eqnarray*}
|(\varphi_1, S^{(0)}_t\varphi_2)|
&=& |(U_\rho^{-1}\varphi_1,(U_\rho S^{(0)}_tU_\rho^{-1})\,U_\rho \varphi_2)|  \\[5pt]
&\leq & e^{-\rho(r-r_1-r_2)} e^{\rho^2 t} \|\varphi_1\|_2 \, \|\varphi_2\|_2 
\;\;\;. 
\end{eqnarray*}
Then setting $\rho= (2t)^{-1} (r-r_1-r_2) > 0$ gives the bounds
\[
|(\varphi_1,S^{(0)}_t\varphi_2)|
\leq e^{-(r-r_1-r_2)^2/(4t)} \, \|\varphi_1\|_2 \, \|\varphi_2\|_2
\]
for all $t > 0$.
Since the estimate is valid for all $r \in \langle0,\infty\rangle$ 
such that $r < d_C(x_1\,;x_2)$ and $r-r_1-r_2 > 0$ one can take the limit 
$r \uparrow d_C(x_1\,;x_2)$ and one obtains 
\[
|(\varphi_1,S^{(0)}_t\varphi_2)|
\leq e^{-(d_C(x_1;x_2)-r_1-r_2)^2/(4t)} \, \|\varphi_1\|_2 \, \|\varphi_2\|_2
\]
for all $t > 0$.\hfill$\Box$

\ruimte

Next we observe that the wave equation associated with $H_0$ has finite speed of
propagation \cite{CGT} \cite{Melr} \cite{Sik}.

\begin{prop} \label{padty3.2}
If  $\varphi_1 \in L_2(B_1)$, $\varphi_2 \in L_2(B_2)$ 
then 
\begin{equation}
(\varphi_1,{\cos(tH_0^{1/2})}\varphi_2)=0
\label{epadty3.1}
\end{equation}
for all $t\leq \tilde d_C(B_1\,;B_2)$.
\end{prop}
\proof\
This  is in fact a corollary of Proposition~\ref{padty311} since the off-diagonal bounds
are equivalent to the finite speed of  propagation by 
the reasoning of \cite{Sik3}.
The principal idea is the following.

\begin{lemma} \label{ladty3.1}
Let $H$ be a positive self-adjoint operator on the Hilbert space $\cal H$.
Fix $\varphi,\psi\in{\cal H}$ and $r\in\langle0,\infty\rangle$.
The following conditions are equivalent.
\begin{tabel}
\item\label{ladty3-10}
$|(\psi,e^{-tH}\varphi)|
\leq e^{-r^2/(4t)}\,\|\psi\|\,\|\varphi\|$
for all $t > 0$.
\item\label{ladty3-1}
There is an $a\geq1$ such that 
$|(\psi,e^{-tH}\varphi)|\leq a\,e^{-r^2/(4t)}\,\|\psi\|\,\|\varphi\|$
for all   $t > 0$.
\item\label{ladty3-2}
$(\psi,{\cos(tH^{1/2})}\varphi)=0$ for all $t\leq r$.
\end{tabel}
\end{lemma}
\proof\
Clearly Condition~\ref{ladty3-10} implies Condition~\ref{ladty3-1}.
Next assume  Condition~\ref{ladty3-1}.
Let $\Ci_+=\{z\in\Ci:\RRe z\geq0,z\neq0\}$ and set $S_z=e^{-zH}$ for all $z \in \Ci_+$.
Define $u \colon \Ci_+ \to \Ri$ by
\[
u(z)=(\psi,S_{z^{-1}}\varphi)
\;\;\; .  \]
Then $u$ is  continuous,
bounded on  its domain of definition and an analytic function on 
$ \{ z \in \Ci : \RRe z > 0 \} $.
Then
\[
\sup_{t > 0}|e^{r^2t/4}u(t)|\leq 
a\,\|\psi\|\,\|\varphi\|
\;\;\;.
\]
Moreover, it follows from positivity and self-adjointness of $H$ that
\[
\sup_{z \in i \Ri \backslash \{ 0 \}  }|e^{r^2z/4}u(z)|
\leq \|\psi\|\,\|\varphi\|
\;\;\;.
\]
 Hence, by the Phragm\'en-Lindel\"of theorem for a quadrant
(see \cite{Mark} vol.~II, Theorem~7.5,  or \cite{SW}, Lemma~4.2, 
or \cite{GS}, Section~IV.7.2), 
\begin{equation}
|e^{r^2z/4}u(z)|  \leq  
a\,\|\psi\|\,\|\varphi\|
\;\;\;.\label{adt1.20}
\end{equation}
Consequently
\begin{equation}
|u(z)|\leq a\, e^{-r^2(\RRe z)/4} 
\|\psi\|\,\|\varphi\|
\label{eod13}
\end{equation}
for all $z$ such that $\RRe z> 0$.
Now
\begin{equation}
(\psi,S_t\varphi)
= (\pi t)^{-1/2}\int^\infty_0ds\,
e^{-s^2/(4t)}(\psi,\cos(sH^{1/2})\varphi)
\;\;\;.
\label{adt1.2}
\end{equation}
Hence changing variables one finds
\[
t^{-1/2}u(4t)=\int^\infty_0ds\,e^{-st} \, w(s)
\]
with $w(s)=(\pi s)^{-1/2}(\psi,\cos(sH)^{1/2}\varphi)$.
Therefore $z\mapsto z^{-1/2}u(4z)$ is the Fourier--Laplace transform of the 
function $s \mapsto w(s)$.
Then it follows from the bounds (\ref{eod13}) and the Paley--Wiener theorem
(\cite{Hor4}, Theorem~7.4.3) that $w$ is supported in the half-line 
$[r^2 ,\infty\rangle$.
Hence Condition~\ref{ladty3-2} is valid.

Now assume Condition~\ref{ladty3-2}.
Then  the integral relation (\ref{adt1.2}) gives
\begin{eqnarray}
|(\psi,S_t\varphi)|
&\leq& (\pi t)^{-1/2}\int^\infty_0ds\,
e^{-s^2/(4t)}|(\psi,\cos(sH^{1/2})\varphi)|\nonumber\\[5pt]
&\leq&(\pi t)^{-1/2}\int^\infty_rds\,
e^{-s^2/(4t)}
\|\psi\|\,\|\varphi\|
\leq e^{-r^2/(4t)}\|\psi\|\,\|\varphi\|
\label{adt1.21}
\end{eqnarray}
and so Condition~\ref{ladty3-10} is valid.
\hfill$\Box$

\ruimte 

The statement of Proposition~\ref{padty3.2} follows immediately from 
Proposition~\ref{padty311} and Lemma \ref{ladty3.1} if $d_C(x_1\,;x_2)<\infty$
and by taking a limit $r\to\infty$ if the distance is infinite.
\hfill$\Box$

\ruimte

The next  lemma is a  simple
consequence of Proposition~\ref{padty3.2}.

\begin{lemma}\label{lod10}
If $\Psi$ is an  even bounded Borel function with 
$\supp \Psi\subseteq [-1,1]$ then
\[
(\varphi_1, {\widehat \Psi(rH_0^{1/2})} \varphi_2)=0
\]
for all $r \leq \tilde d_C(B_1\,;B_2)$, $\varphi_1 \in L_2(B_1)$ and $\varphi_2 \in L_2(B_2)$,
where \raisebox{0pt}[0pt][0pt]{$\widehat \Psi$} denotes the 
Fourier transform of $\Psi$.
\end{lemma}
\proof\
Since $\Psi$ is even,
\[
\widehat \Psi(rH_0^{1/2}) 
=(2\pi)^{-1}\int_{\Ri}dt\, \Psi(t) \cos(rtH_0^{1/2})\;\;\;.
\]
But  $\supp \Psi \subseteq [-1,1]$ 
and the statement of the lemma follows immediately from Proposition~\ref{padty3.2}.
\hfill$\Box$

\ruimte

The $L_2$ off-diagonal bounds can be extended to more general sets
than balls  for strongly elliptic operators or
operators with continuous coefficients but this is not strictly relevant to the sequel.
One can also  derive off-diagonal bounds for general sets without assuming strong ellipticity
or continuity of the coefficients if one uses the Euclidean distance.
Then, however, the Gaussian factor changes.

\begin{lemma}\label{ladt212}
If  $\|C\|\neq0$  then 
\[
|(\varphi_1,S^{(0)}_t\varphi_2)|
\leq e^{- d_e(V_1;V_2)^2/(4\|C\|t)}\, \|\varphi_1\|_2 \, \|\varphi_2\|_2
\]
for all non-empty measurable $V_1,V_2\subset\Ri^d$ and all
$\varphi_1\in L_2(V_1)$, $\varphi_2\in L_2(V_2)$,  $t > 0$, 
where $d_e(V_1\,;V_2)$ is
the Euclidean distance between $V_1$ and $V_2$.
\end{lemma}
\proof\
Set $N = d_e(V_1\,;V_2) + 1$ and define $\psi \colon \Ri^d \to \Ri$ by
$\psi(x) = d_e(x\,;V_2) \wedge N$.
Then $\psi$ is bounded and $|\psi(x) - \psi(y)| \leq |x-y|$ for 
all $x,y \in \Ri^d$.
Therefore $\psi$ is partial differentiable, in the $L_\infty$ sense, and 
$\sum_{i,j=1}^d (\partial_i \psi) \, c_{ij} \, (\partial_j \psi) 
\leq \|C\|\sum^d_{i=1}|\partial_i\psi|^2\leq \|C\|$ almost everywhere.

Now for all $\rho \in \Ri$ define the bounded multiplication operator
$U_\rho$ by $U_\rho \varphi = e^{-\rho \psi} \varphi$.
Then one computes as in the proof of Proposition~\ref{padty311} that 
$\|U_\rho S^{(0)}_t U_\rho^{-1}\|_{2 \to 2} \leq e^{\|C\| \rho^2 t}$
for all $\rho \in \Ri$ and $t > 0$.
Next,
\[
|(\varphi_1, S^{(0)}_t \varphi_2)|
= |(U_\rho^{-1} \varphi_1, (U_\rho S^{(0)}_t U_\rho^{-1}) U_\rho \varphi_2)|
\leq e^{\|C\| \rho^2 t} \|U_\rho^{-1} \varphi_1\|_2 \, \|U_\rho \varphi_2\|_2
\;\;\; .  \]
But $U_\rho \varphi_2 = \varphi_2$ and if $\rho \leq 0$
then 
$\|U_\rho^{-1} \varphi_1\|_2 
      \leq e^{\rho (d_e(V_1;V_2) \wedge N)} \|\varphi_1\|_2 
      = e^{\rho d_e(V_1;V_2)}\|\varphi_1\|_2$.
Hence choosing $\rho = - (2t)^{-1} d_e(V_1;V_2)$ establishes the lemma.\hfill$\Box$

\ruimte

The last lemma allows one to prove that $S^{(0)}$ conserves probability, i.e.,
the Dirichlet form $h_0$ is conservative in the terminology of \cite{FOT}, page~49.

\begin{prop}\label{pod1}
The extension of $S^{(0)}$ to $L_\infty(\Ri^d)$ satisfies
\[
S^{(0)}_t\one=\one
\]
for all $t>0$.
\end{prop}
\proof\
We may assume that $\|C\| \neq 0$.
Let $\varphi \in C_c(\Ri^d)$.
Let $R > 0$  and suppose that $\supp \varphi \subset B_R$ where $B_R$ is the 
(Euclidean) ball of 
radius $R$ centred at the origin.
Let $\chi_R$ be a positive $C^\infty$-function with $\chi_R(x)=1$ if $|x|\leq 2R$
and $\chi_R(x)=0$ if $|x|\geq 3R$.
Then with $(\cdot,\cdot)$ denoting as usual the pairing between $L_p$ and $L_q$ one has
\begin{eqnarray*}
|(S^{(0)}_t\one,\varphi)-(\one,\varphi)|
&\leq&
|(S^{(0)}_t(\one-\chi_R),\varphi)|
+|(S^{(0)}_t\chi_R,\varphi)
-(S^{(\varepsilon)}_t\chi_R,\varphi)|
\nonumber\\[5pt]
&&\hspace{2cm}{}+|(S^{(\varepsilon)}_t(\one-\chi_R),\varphi)|
\end{eqnarray*}
 for all $\varepsilon > 0$, where we have used $S^{(\varepsilon)}_t\one=\one$.
The latter equality follows from the strong ellipticity of $H_\varepsilon$ (see, for example,
\cite{ER22}, page~145, proof of Theorem~4.6).
Now Lemma~\ref{ladt212} gives
\begin{eqnarray*}
|(S^{(0)}_t(\one-\chi_R),\varphi)|
& \leq & \sum_{n=2}^\infty 
       |(\one_{B_{(n+1)R} \backslash B_{nR}} (\one-\chi_R),S^{(0)}_t \varphi)|  \\[5pt]
& \leq & \sum_{n=2}^\infty e^{-d_e(B_{(n+1)R} \backslash B_{nR} ;B_R)^2 (4\|C\|t)^{-1}}
     \|\one_{B_{(n+1)R} \backslash B_{nR}} \|_2
     \|\varphi\|_2  \\[5pt]
& \leq & \sum_{n=2}^\infty e^{-(n-1)^2R^2(4\|C\|t)^{-1}}
     |B_{(n+1)R}|^{1/2}
     \|\varphi\|_2  \\[5pt]
& \leq & \sum_{n=2}^\infty (n+1)^{d/2}e^{-(n-1)^2R^2(4\|C\|t)^{-1}}|B_R|^{1/2}
     \|\varphi\|_2  
\;\;\; .
\end{eqnarray*}
Similarly, 
\[
|(S^{(\varepsilon)}_t(\one-\chi_R),\varphi)|
\leq \sum_{n=2}^\infty (n+1)^{d/2}e^{-(n-1)^2R^2(4\|C\|t)^{-1}}|B_R|^{1/2}
     \|\varphi\|_2  
\]
for all $\varepsilon > 0$.
Hence
\begin{eqnarray*}
|(S^{(0)}_t\one,\varphi)-(\one,\varphi)|
& \leq & |(S^{(0)}_t\chi_R,\varphi)
-(S^{(\varepsilon)}_t\chi_R,\varphi)|  \\*[5pt]
& & \hspace{20mm} {}
+ 2  \sum_{n=2}^\infty (n+1)^{d/2}e^{-(n-1)^2R^2(4\|C\|t)^{-1}}|B_R|^{1/2}
     \|\varphi\|_2 
\end{eqnarray*}
for all $\varepsilon > 0$ and $R > 0$ such that $\supp \varphi \subset B_R$.
Since $S^{(\varepsilon)}_t$ converges strongly to $S^{(0)}_t$ as $\varepsilon\to0$
the desired result follows by taking successive limits $\varepsilon\to0$ 
and 
$R\to\infty$.\hfill$\Box$

\ruimte

The proof of the conservation property of Proposition~\ref{pod1} is partially based on 
the observation that it is valid for semigroups generated by strongly elliptic operators.
But in the latter case the statement of the proposition 
extends to a wider class of functions by general functional  analysis.

\begin{cor}\label{cadt213}
Let $H$ be a strongly elliptic operator in divergence form with real measurable 
coefficients and $\Phi$  a function which is bounded and holomorphic in a strip 
$\{ z\in\Ci : |\IIm z|<2\delta \} $ for some $\delta > 0$.
Then $\Phi(H^{1/2})$ extends to a bounded operator on $L_\infty$ and 
the extension, still denoted by $\Phi(H^{1/2})$, satisfies
\begin{equation}
\Phi(H^{1/2})\one=\Phi(0)\one\;\;\;.
\label{eadt2.99}
\end{equation}
\end{cor}
\proof\
If $S$ denotes the self-adjoint semigroup generated by $H$ on $L_2(\Ri^d)$ 
then $S$ extends to a positive contraction semigroup on the spaces $L_p(\Ri^d)$ 
for all $p\in[1,\infty]$, and $S_t\one=\one$.
Moreover the semigroup $S$ on $L_2(\Ri^d)$ is holomorphic in 
the open right half-plane.
But its kernel satisfies Gaussian bounds and these extend to the open right half-plane 
(see \cite{Dav2} especially Theorem~3.4.8, or \cite{ER18} Theorem~1.1). 
Therefore the extension of $S$ to the $L_p$-spaces is also holomorphic 
in the open right half-plane.

Next observe that the Poisson semigroup  $P$  generated by $H^{1/2}$ is given by
\[
P_t= (4\pi )^{-1/2}\int^\infty_0ds\,t \, s^{-3/2}e^{-t^2/(4s)}S_s
\;\;\;.
\]
Therefore $P_t$ maps $L_\infty$ into $L_\infty$ and $P_t \one = \one$ for all $t > 0$.
Moreover it readily follows that $P$ is holomorphic in the open right half-plane
on each of the $L_p$-spaces.
Then $(\lambda I + H^{1/2})^{-1}$ maps $L_\infty$ into $L_\infty$ 
for all $\lambda \in \Ci$ with $\RRe \lambda > 0$
and $(\lambda I + H^{1/2})^{-1} \one = \lambda^{-1} \one$.
In addition $z \mapsto (z I + H^{1/2})^{-1}$ is analytic on $\Ci \backslash \langle-\infty,0]$.
Hence $(\lambda I + H^{1/2})^{-1} \one = \lambda^{-1} \one$ for all
$\lambda \in \Ci \backslash \langle-\infty,0]$.

If $\Gamma$ is a contour in the complex plane from $\infty + i \delta$ to 
$\infty - i \delta$
contained in the set $\{z\in\Ci: |\IIm z|<2\delta \} \backslash [0,\infty\rangle$
then
\[
(2\pi i)^{-1} \int_\Gamma dz \, \Phi(z) \, (z I - H^{1/2})^{-1}
\]
is an operator which maps $L_\infty$ into $L_\infty$ and $L_2$ into $L_2$.
Moreover, it equals the operator $\Phi(H^{1/2})$ defined on $L_2$ by spectral theory.
Hence $\Phi(H^{1/2})$ extends to a bounded operator on $L_\infty$.
In addition 
\[
\Phi(H^{1/2}) \, \one
= (2\pi i)^{-1} \int_\Gamma dz \,\Phi(z) \, z^{-1} \, \one
= \Phi(0) \, \one
\]
as desired.\hfill$\Box$

\section{Subelliptic heat kernel estimates}\label{Sadty4}

The $L_2$ off-diagonal bounds  derived in Proposition~\ref{padty311} are valid for all
second-order elliptic operators in divergence form.
Next we examine pointwise Gaussian upper  bounds on the distribution kernel $K^{(0)}$ of the 
viscosity semigroup $S^{(0)}$ but for these we require a  subellipticity
estimate.
We define $H_0$ to be subelliptic (of order $\gamma$) if there exist $\mu>0$, 
$\nu\geq0$ and $\gamma\in \langle0,1]$
such that 
\begin{equation}
H_0\geq \mu \,\Delta^\gamma-\nu\, I
\;\;\;.
\label{eadty4.28}
\end{equation}
in the sense of quadratic forms.
This is equivalent to requiring that the Sobolev inequalities
\[
\|H_0^{1/2}\varphi\|_2^2+\nu\,\|\varphi\|_2^2\geq\mu\,\|\Delta^{\gamma/2}\varphi\|_2^2
\]
are satisfied for all $\varphi\in D(H_0^{1/2})=D(h_0)$.
Note that if $\gamma=1$ then  one may choose $\nu=0$ and (\ref{eadty4.28}) 
reduces to the strong ellipticity
condition $H_0\geq\mu\,\Delta$ (see Remark~\ref{radty2.1}).
Note further that the order is not uniquely defined. 
If (\ref{eadty4.28}) is satisfied for one value of $\gamma$ it is also satisfied for all smaller 
$\gamma\in\langle0,1\rangle$ since one has inequalities $\Delta^\alpha\leq a\,(\Delta^\beta+I)$ 
for $\beta>\alpha$.

In the sequel when we assume that $H_0$ is subelliptic then $\mu$, $\nu$ and 
$\gamma$ will always denote the
parameters in the subellipticity condition (\ref{eadty4.28}).

The first result is an estimate for small $t$ which follows by variation of standard arguments.

\begin{prop} \label{padty4.1}
Assume the viscosity operator $H_0$ is subelliptic of order $\gamma$.

There is an $a>0$, depending only on $\gamma$ and $d$, such that
\begin{equation}
\|S^{(0)}_t\|_{p \to q} 
\leq a \, e^\nu\,(\mu (t\wedge 1))^{-d/(2r\gamma)} 
\label{eeo4.5;32}
\end{equation}
for all $t > 0$  and $p,q \in [1,\infty]$ with $p \leq q$ where $r^{-1}=p^{-1}-q^{-1}$.

Moreover, for each $\delta>0$ there exists an  $a>0$, depending only on $\delta$, 
$\gamma$ and $d$, such that
\begin{equation}
K^{(0)}_t(x\,;y)
\leq a \, e^{2\nu}\,(\mu (t\wedge 1))^{-d/(2\gamma)} \,  
e^{-{d_C}(x;y)^2 ((4+\delta)t)^{-1}}
\;\;\;\; (x,y) \mbox{-a.e.}
\label{eeo4.5}
\end{equation}
uniformly for all $t>0$.
\end{prop}
\proof\
The starting point is the fractional Nash inequality 
\begin{equation}
\|\varphi\|_2^{2+4 \gamma/ d}
\leq c_1 \, \|\Delta^{\gamma/2} \varphi\|_2^2 \,\|\varphi\|_1^{4 \gamma/ d}
\label{eeo4.01}
\end{equation}
which is valid for all $\varphi \in L_1(\Ri^d) \cap D(\Delta^{\gamma/2})$ and a $c_1>0$, depending 
only on $\gamma$ and~$d$.
This follows by a slight variation of Nash's original arguments (which he attributes to
Stein, see \cite{Nash} page~935)  for 
strongly elliptic operators (see, for example, \cite{Robm} page~169). 
The principal point is that if $T$ denotes the self-adjoint semigroup generated by $\Delta^\gamma$ 
then the H\"older inequality gives
\[
\|T_t\|_{1\to2}
=\|T_t\|_{2\to\infty}
=\bigg(\int_{\Ri^d}dp\,e^{-2tp^{2\gamma}}\bigg)^{1/2}
= t^{-d/(2\gamma)}\|T_1\|_{2\to\infty}
\]
for all $t>0$.
Therefore
\begin{eqnarray*}
\|\varphi\|_2&\leq&\|(I-T_t)\varphi\|_2+\|T_t\|_{1\to2}\|\varphi\|_1\\[5pt]
&\leq&t\,\|\Delta^{\gamma/2}\varphi\|_2+t^{d/(2\gamma)}\,\|T_1\|_{1\to2}\|\varphi\|_1
\end{eqnarray*}
for all $t>0$.
Optimization with respect to $t$ then yields (\ref{eeo4.01}).

Combination of the Nash inequality (\ref{eeo4.01}) and the subellipticity condition (\ref{eadty4.28})
gives the  Nash inequality
\begin{equation}
\|\varphi\|_2^{2+4 \gamma /d}
\leq c_1 \, \mu^{-1}  \Big(h_0(\varphi)+\nu\, \|\varphi\|_2^2\Big) \|\varphi\|_1^{4 \gamma / d}
\label{eeo4.2}
\end{equation}
for all $\varphi \in L_1(\Ri^d) \cap D(h)$.

Next in order to avoid domain problems
we use the approximants $h_\varepsilon$.
Since  $h_\varepsilon\geq h_0$ there are  inequalities  similar to (\ref{eeo4.2}) with $h_0$
replaced by $h_\varepsilon$.
Then, following Nash (see, for example, \cite{CKS} Theorem~2.1), one obtains 
bounds
\begin{equation}
\|S^{(\varepsilon)}_t\|_{1 \to 2}
=\|S^{(\varepsilon)}_t\|_{2 \to \infty}
\leq  c_2 \, (\mu t)^{- d/(4\gamma)} \, e^{\nu t}
\label{eadty4.73}
\end{equation}
uniform for all $\varepsilon, t>0$, where $c_2=(cd/(4\gamma))^{d/(4\gamma)}$.
Then it follows from the strong convergence that 
\[
\|S^{(0)}_t \varphi\|_2
= \lim_{\varepsilon \downarrow 0} \|S^{(\varepsilon)}_t \varphi\|_2
\leq c_2 \, (\mu t)^{- d/(4\gamma)} \, e^{\nu t} \|\varphi\|_1
\]
for all $t > 0$ and $\varphi \in L_1 \cap L_2$.
Since $S^{(0)}$ is a contraction semigroup it further follows that 
$\|S^{(0)}_t\varphi\|_2\leq \|S^{(0)}_1\varphi\|_2$ for all $t\geq 1$  and $\varphi \in L_1 \cap L_2$.
Hence 
\[
\|S^{(0)}_t\|_{1 \to 2} \leq c_2 \,e^\nu\, (\mu (t\wedge 1))^{- d/(4\gamma)}
\]
for all $t>0$.
Then the bounds (\ref{eeo4.5;32}) follow with the aid of the contractivity of $S^{(0)}$ 
 by interpolation.

One also has  bounds analogous to (\ref{eeo4.5;32}) for the approximants $S^{(\varepsilon)}$ and so  
\[
\|K^{(\varepsilon)}_t\|_\infty
=\|S^{(\varepsilon)}_t\|_{1 \to \infty}
\leq \|S^{(\varepsilon)}_{t/2}\|_{2\to\infty}^2
\leq  c_3 \, (\mu (t\wedge 1))^{- d/(2\gamma)} 
\]
uniformly for all $\varepsilon, t>0$ with $K^{(\varepsilon)}$ the kernel of the  
approximating semigroups $S^{(\varepsilon)}$ and $c_3=c_2^2\,e^{2\nu}\,2^{d/(2\gamma)}$.
Now one can extend  these latter bounds
to the  Davies' perturbation of $S^{(\varepsilon)}$ by Davies' method \cite{Dav7}
as elaborated by Fabes and Stroock  \cite{FaS} to obtain the Gaussian bounds
\begin{equation}
K^{(\varepsilon)}_t(x\,;y)
\leq a \,  e^{2\nu}\,(\mu ( t\wedge 1))^{-d/(2\gamma)}
\,e^{-d_{C_\varepsilon}(x;y)^2 ((4+\delta)t)^{-1}}
\label{eeo4.5.1}
\end{equation}
with $a$ independent of $\varepsilon$,   $\mu$ and $\nu$.

Next one has the following convergence result for the kernels.

\begin{lemma}\label{leo4.1}
Assume   the viscosity operator $H_0$ is subelliptic.
Then the kernels $K^{(\varepsilon)}_t$ 
converge in the weak$^*$ sense on $L_\infty(\Ri^d\times\Ri^d)$,
as $\varepsilon\to0$, to the kernel $K^{(0)}$.
\end{lemma}
\proof\
If $\varepsilon\in\langle0,1]$ then $\|C_\varepsilon\|\leq\|C\|+1$.
Therefore $d_{C_\varepsilon}(x\,;y)\geq (1+\|C\|)^{-1/2}|x-y|$ and it follows from (\ref{eeo4.5.1}) that
there are $a,b > 0$, depending only on $\mu$, $\nu$, $\gamma$ and $\|C\|$, such that
\begin{equation}
K^{(\varepsilon)}_t(x\,;y)
\leq a \, (t\wedge 1)^{-d/(2\gamma)} \,  
e^{-b|x-y|^2 t^{-1}}
\label{eeo4.5.2}
\end{equation}
uniformly for all $t>0$, $x,y \in \Ri^d$ and $\varepsilon\in \langle0,1]$.
The convergence of the  $K^{(\varepsilon)}$ follows from these uniform upper bounds 
and the $L_2$-convergence of $S^{(\varepsilon)}$ to
$S^{(0)}$  (see, for example, \cite{ER15} proof of Proposition~2.2).
\hfill$\Box$

\ruimte

Moreover one has convergence of the distances.

\begin{lemma}\label{ladty4.2}
$\displaystyle \lim_{\varepsilon \to 0} d_{C_\varepsilon}(x\,;y)
= \sup_{\varepsilon > 0} d_{C_\varepsilon}(x\,;y)
= d_C(x\,;y)$  
for all $x,y \in \Ri^d$.
\end{lemma}
\proof\
If $\varepsilon_1 \leq \varepsilon_2$ then 
$d_{C_{\varepsilon_2}} \leq d_{C_{\varepsilon_1}} \leq d_C$ by the definition
of the distances.
So 
\begin{equation}
\lim_{\varepsilon \to 0} d_{C_\varepsilon}(x\,;y)
= \sup_{\varepsilon > 0} d_{C_\varepsilon}(x\,;y)
\leq d_C(x\,;y)
\label{epdist8.3;1}
\end{equation}
for all $x,y \in \Ri^d$.
Fix $x,y \in \Ri^d$.
Let $\psi \in C_c^\infty(\Ri^d)$ and suppose 
$\sum_{i,j=1}^d (\partial_i \psi) \, c_{ij} \, (\partial_j \psi) \leq 1$
almost everywhere.
Then 
$\sum_{i,j=1}^d (\partial_i \psi) \, (c_{ij} + \varepsilon \delta_{ij}) \, (\partial_j \psi) 
    \leq 1 + \varepsilon M$,
with $M = \sum_{i=1}^d \|\partial_i \psi\|_\infty^2 $,
 almost everywhere for all $\varepsilon > 0$.
If $\psi_\varepsilon = (1 + \varepsilon M)^{-1/2} \psi$ then 
$\psi_\varepsilon \in C_c^\infty(\Ri^d)$ and 
$\sum_{i,j=1}^d (\partial_i \psi_\varepsilon) \, (c_{ij} + \varepsilon \delta_{ij}) \, (\partial_j \psi_\varepsilon) 
    \leq 1$ almost everywhere.
So 
\[
d_{C_\varepsilon}(x\,;y) 
\geq |\psi_\varepsilon(x) - \psi_\varepsilon(y)|
= (1 + \varepsilon M)^{-1/2} \, |\psi(x) - \psi(y)|
\]
for all $\varepsilon > 0$.
Now take the limit $\varepsilon \to 0$.
Then 
\[
\lim_{\varepsilon \to 0} d_{C_\varepsilon}(x\,;y)
\geq |\psi(x) - \psi(y)|
\;\;\; .  \]
But this implies that 
$\lim_{\varepsilon \to 0} d_{C_\varepsilon}(x\,;y) \geq d_C(x\,;y)$
and the lemma follows.\hfill$\Box$

\ruimte

The statement of Proposition~\ref{padty4.1} concerning the kernel follows from 
(\ref{eeo4.5.1}) in the limit $\varepsilon\to0$
as a consequence of Lemmas~\ref{leo4.1} and \ref{ladty4.2}.\hfill$\Box$

\ruimte

There are a number of  alternative ways of passing from the semigroup estimates to the pointwise estimates
on the kernel.
Theorem~4 of \cite{Sik3} is based on an argument which exploits the finite speed of propagation and 
which is applicable in the current context.

\smallskip

Note that if for some $x,y\in\Ri^d$ one has ${d_C}(x;y)=\infty$ 
and $K^{(0)}_t$ is continuous at $(x,y)$ then $K^{(0)}_t(x\,;y) = 0$.
Further the foregoing  arguments  give a  bound on the
 kernel which does not decrease with $t$. 
We will see in Section~\ref{Sadty5} that this is the best one can hope
for unless one  has more information such as continuity or strict positivity
of the kernel.

\smallskip

The Gaussian upper bounds in fact give information on lower bounds by a variation of standard arguments
for strongly elliptic operators.

\begin{cor}\label{ctdy3.3}
Assume that the viscosity operator $H_0$ is subelliptic.
Let $r,t>0$. 
Then there is an  $a'>0$ such that
\begin{equation}
(\varphi,S^{(0)}_t\varphi)\geq a'\,\|\varphi\|_1^2
\label{etdy3.20}
\end{equation}
for all positive $\varphi\in L_1(\Ri^d)\cap L_2(\Ri^d)$ with $\diam(\supp \varphi)\leq r$,
where the diameter is with respect to the Euclidean distance.
Hence if $K^{(0)}_t$ is continuous at $(x,x)\in\Ri^d\times \Ri^d$ then
\[
K^{(0)}_t(x\,;x)\geq a'
\;\;\;.
\]
The value of $a'$ depends on $H_0$ only through the parameters $\mu$, $\nu$, $\gamma$ and~$\|C\|$.
\end{cor}
\proof\
Since $S^{(0)}$ is self-adjoint it follows that 
$(\varphi,S^{(0)}_t\varphi)=\|S^{(0)}_{t/2}\varphi\|_2^2\geq 0$.
 Hence
\begin{equation}
|(\varphi, S^{(0)}_t\psi)|^2
\leq (\varphi, S^{(0)}_t\varphi)\,(\psi, S^{(0)}_t\psi)
\label{etdy3.20.1}
\end{equation}
for all $\varphi,\psi\in L_2(\Ri^d)$.
Next  let $x_0 \in \Ri^d$ and let
$\varphi$ be a positive integrable function with support in the Euclidean ball 
$B_e(x_0\,;r)=\{y\in \Ri^d: |y-x_0|<r\}$.
Further let $R > 2r$ and let $\psi$ be the characteristic function of the ball $B_e(x_0\,;R)$.
We  evaluate (\ref{etdy3.20.1}) with this choice of $\varphi$ and $\psi$.

First one has
\[
(\psi, S^{(0)}_t\psi)\leq \|\psi\|_2^2=V_e(R)
\]
where $V_e(R)$ is the volume of $B_e(x_0\,;R)$.

Secondly,   $S^{(0)}_t\one=\one$.
Hence
\[
(\varphi, S^{(0)}_t\psi)
=(\varphi,\one)-(\varphi,S^{(0)}_t(\one-\psi))
\geq \|\varphi\|_1 \Big(1- \sup_{x \in B_e(x_0;r)} \int_{ \{ y : |y-x_0|\geq R \} } dy\, 
                                         K^{(0)}_t(x\,;y)\Big)
\;\;\;.
\]
Then since $K^{(\varepsilon)}_t$ satisfies  the bounds (\ref{eeo4.5.2}) 
it follows from Lemma~\ref{leo4.1} that there are 
$a,b > 0$, depending only on $\mu$, $\nu$, $\gamma$ and $\|C\|$, such that
\begin{equation}
K^{(0)}_t(x\,;y)
\leq a \, ( t\wedge 1)^{-d/(2\gamma)} \,  e^{-b|x-y|^2 t^{-1}}
\label{ectdy3.3;1}
\end{equation}
uniformly for all $t>0$, $x,y \in \Ri^d$.
Hence one can choose $R$
sufficiently large that
\[
(\varphi, S^{(0)}_t\psi)
\geq 2^{-1} \|\varphi\|_1
\;\;\; .  \]
Thirdly,  substituting these last two estimates in (\ref{etdy3.20.1}) one deduces that 
\[
(\varphi,S^{(0)}_t\varphi)
\geq (4V_e(R))^{-1} \, \|\varphi\|_1^2 
\;\;\; .  \]
It follows  that (\ref{etdy3.20}) is valid with $a'= (4V_e(R))^{-1}$.
The value of $R$ is dictated by the Gaussian bounds (\ref{ectdy3.3;1}) and hence depends on
$H_0$ only through the parameters $\mu$, $\nu$, $\gamma$ and~$\|C\|$.

Finally suppose $K^{(0)}_t$ is continuous at a diagonal point, which we may take to be
$(0,0)$. 
Then for $\lambda > 0$ replace $\varphi$ in (\ref{etdy3.20}) by $\varphi_\lambda$ where 
$\varphi_\lambda(x)=\lambda^{-d}\varphi(\lambda^{-1}x)$.
It follows that $\|\varphi_\lambda\|_1=\|\varphi\|_1$.
Moreover,
\[
\lim_{\lambda\to0}(\varphi_\lambda,S^{(0)}_t\varphi_\lambda)=
\lim_{\lambda\to0}\int_{\Ri^d}dx\int_{\Ri^d}dy\,\varphi(x)\,\varphi(y)\,K^{(0)}_t(\lambda x\,;\lambda y)
=\|\varphi\|_1^2\,K^{(0)}_t(0\,;0)
\;\;\;.
\]
Therefore $K^{(0)}_t(0\,;0)\geq a'$.\hfill$\Box$ 

\begin{remarkn}\label{reo4.1}
If the kernel $K^{(0)}_t$ is a continuous function on $\Ri^d\times \Ri^d$ it follows
from the corollary that it is strictly positive on the diagonal, i.e., 
\[
\inf_{x\in\Ri^d}K^{(0)}_t(x\,;x)\geq a'>0
\;\;\;.
\]
If, however, $K^{(0)}_t$ is uniformly continuous then one has a stronger
off-diagonal property.
Explicitly, if 
\[
\lim_{x\to0}\|L(x)K^{(0)}_t-K^{(0)}_t\|_\infty=0
\]
where   $(L(z)K^{(0)}_t)(x\,;y)=K^{(0)}_t(x-z\,;y)$
then it follows from Corollary~\ref{ctdy3.3} that there are $a',r>0$ such that 
\[
K^{(0)}_t(x\,;y)\geq a'>0
\]
for all $x, y\in\Ri^d$ with $|x-y|<r$.
Uniform continuity of the kernel in the first variable is of course equivalent to uniform
continuity in the second variable, by symmetry, and separate uniform continuity is
equivalent to joint uniform continuity.
\end{remarkn}

If the kernel is uniformly continuous then subellipticity implies large time Gaussian bounds
of a different character geometric character to the small time bounds of Proposition~\ref{padty4.1}.
The uniform continuity implies that the kernel decays as $t\to\infty$ with the rate of decay dictated
by the dimension $d$ independent of the order of subellipticity.

\begin{thm}\label{tadty411}
Let $H_0$ be the viscosity operator with coefficients $C=(c_{ij})$
and $K^{(0)}$ the distribution kernel of the contraction semigroup $S^{(0)}$
generated by $H_0$.
Assume 
\begin{enumerate}
\item\label{padty4.2-1}
$H_0$ is subelliptic.
\item \label{padty4.2-2}
There are $a, r>0$ such that $K^{(0)}_1(x\,;y)\geq a$
for almost every $(x,y)\in \Ri^d \times \Ri^d$ with $|x-y|\leq r$.
\end{enumerate}
Then for all $\delta>0$  there exists an  $a'>0$ such that
\begin{equation}
K^{(0)}_t(x\,;y)
\leq a' \, t^{-d/2} 
e^{-{d_C}(x;y)^2 ((4+\delta)t)^{-1}}
\;\;\;, 
\label{etdy4.111}
\end{equation}
$(x,y)$ almost everywhere, for all  $t\geq1$.

Moreover, for all $R>0$  there is an  $a''>0$ such that
\begin{equation}
(\varphi,S^{(0)}_t\varphi)\geq a''\,t^{-d/2}\,\|\varphi\|_1^2
\label{etdy3.202}
\end{equation}
for all $t \geq 1$ and positive $\varphi\in L_1(\Ri^d)\cap L_2(\Ri^d)$ with 
$\diam(\supp \varphi)\leq R$.
Hence if $K^{(0)}_t$ is continuous at $(x,x)\in\Ri^d\times \Ri^d$ then
\[
K^{(0)}_t(x\,;x)\geq a''\,t^{-d/2}
\]
for all $t\geq1$.
\end{thm}
\begin{remarkn}\label{reo4.2}
The local lower bounds of Condition~2 follow from 
 Condition~1 if $K^{(0)}_t$ is uniformly continuous. 
This is a consequence of Remark~\ref{reo4.1}.
Nevertheless we show in Section~\ref{Sadty5} that subellipticity does not necessarily
imply uniform continuity nor does it imply strict positivity of the kernel.
\end{remarkn}
\proof\
The proof of the upper bounds is again  based on Nash's original arguments as elaborated
by Carlen, Kusuoka and Stroock \cite{CKS}.     
In particular the following lemma is a version of an argument in Section~4 of \cite{CKS}. 

\begin{lemma}\label{ladty4.3}
Assume that the distribution  kernel $K^{(0)}$  satisfies the local lower
bounds of  Condition~$\ref{padty4.2-2}$ of Theorem~$\ref{tadty411}$.
Then there exists a $\rho>0$ such that
\[
H_0\geq\rho\,\Delta(I+\Delta)^{-1}
\;\;\;.
\]
\end{lemma}
\proof\
Using  spectral theory and  the conservation property of Proposition~\ref{pod1}
one has 
\begin{eqnarray*}
h_0(\varphi)
&\geq &t^{-1}(\varphi, (I-S^{(0)}_t)\varphi)\\[5pt]
& = & (2t)^{-1}\Big( (S^{(0)}_t\one,|\varphi|^2)+(|\varphi|^2,S^{(0)}_t\one)
    - (\varphi, S^{(0)}_t\varphi) - (S^{(0)}_t\varphi,\varphi) \Big) \\[5pt]
& = &(2t)^{-1}\int_{\Ri^d}dx\int_{\Ri^d}dy\,K^{(0)}_t(x\,;y)|\varphi(x)-\varphi(y)|^2
\end{eqnarray*}
for all $\varphi\in D(h_0)$ and $t>0$.
Next choose a smooth positive function $\rho$ with support in $\langle-r,r\rangle$
such that $\rho\leq a$ and  $\rho=a$ if $|x|\leq r/2$.
Then it follows by assumption that
\[
K^{(0)}_1(x\,;y)\geq \,\rho(|x-y|^2 )
\]
for all $x,y\in\Ri^d$.
Combining these inequalities one finds
\[
h_0(\varphi)\geq
\int_{\Ri^d}d\xi\, |{\widehat\varphi}(\xi)|^2 \, \int_{\Ri^d}dx \, \rho(|x|^2)\,
(1-\cos\xi.x)
\;\;\;.
\]
But by the choice of $\rho$ one can find a $\sigma>0$ such that
\[
\int_{\Ri^d}dx\,\rho(|x|^2)\,
(1-\cos\xi.x)\geq\sigma\,(|\xi|^2\wedge1)\geq \sigma\, |\xi|^2(1+|\xi|^2)^{-1}
\]
for all $\xi\in \Ri^d$.
Therefore $h_0(\varphi)\geq \sigma\,l((I+\Delta)^{-1/2}\varphi)$ for all $\varphi\in D(h_0)$.
\hfill$\Box$

\ruimte

Lemma~\ref{ladty4.3} implies that there exists a $\sigma>0$ such that
\begin{equation}
h_0(\varphi)\geq
\sigma \int_{\Ri^d}d\xi\, |{\widehat\varphi}(\xi)|^2\,\left(|\xi|^2\wedge1\right)
\label{eeo4.100}
\end{equation}
for all $\varphi\in D(h)$.
Assume  that $\varphi\in D(h)\cap  L_1$. 
It then   follows by Fourier transformation,
as in the proof of Corollary~4.9 in \cite{CKS},  that
\begin{eqnarray*}
\|\varphi\|_2^2
&=&\int_{ \{ \xi : |\xi|\leq R \} }d\xi\,|\widehat\varphi(\xi)|^2+
\int_{ \{ \xi : |\xi|\geq R \} }d\xi\,|\widehat\varphi(\xi)|^2\\[5pt]
&\leq&c\,R^d\,\|\varphi\|_1^2+\int_{ \{ \xi : |\xi|\geq R \} }d\xi\,(R^{-2}|\xi|^2\wedge 1)\,
|\widehat\varphi(\xi)|^2\\[5pt]
&\leq &c\,R^d\,\|\varphi\|_1^2 + R^{-2}\int_{\Ri^d}d\xi \, (|\xi|^2 \wedge 1) \, |\widehat\varphi(\xi)|^2  \\[5pt]
& \leq & c\,R^d\,\|\varphi\|_1^2+ R^{-2} \sigma^{-1} h_\varepsilon(\varphi)\hspace{1cm}
\end{eqnarray*}
for all $R\in\langle0,1]$ and $\varepsilon > 0$ where the last inequality uses (\ref{eeo4.100})
and $c$ is the volume of the Euclidean unit ball in $\Ri^d$. 
Then  the Nash inequality
\begin{equation}
\|\varphi\|_2^{2+4/d}\leq c'\,h_\varepsilon(\varphi)\|\varphi\|^{4/d}_1
\label{etdy5.31}
\end{equation}
follows for all $\varphi\in D(h_\varepsilon)\cap L_1$ with 
$h_\varepsilon(\varphi)\leq \|\varphi\|^2_1$ by setting 
$R=(h_\varepsilon(\varphi)/\|\varphi\|^2_1)^{1/(d+2)}$.
The inequality is uniform for $\varepsilon\in\langle0,1]$.
Note that (\ref{etdy5.31}) is analogous to the earlier Nash inequality (\ref{eeo4.2}) but
with $\gamma=1$ and $\nu=0$.
In addition there is the important restriction $h_\varepsilon(\varphi)\leq \|\varphi\|^2_1$.

Next it follows from the contractivity of $S^{(\varepsilon)}$ on $L_1$ that 
\begin{equation}
\|S^{(\varepsilon)}_t\|_{1\to\infty}\leq \|S^{(\varepsilon)}_1\|_{1\to\infty}
\label{etdy5.4}
\end{equation}
 for all $t\geq 1$.
In particular $t\mapsto\|S^{(\varepsilon)}_t\|_{1\to\infty}$ is uniformly bounded for $t\geq 1$.
The conditions (\ref{etdy5.31}) and (\ref{etdy5.4}) correspond to the assumptions of Theorem~2.9 of \cite{CKS}.
Therefore the theorem establishes the 
 large time estimates
\begin{equation}
\|S^{(\varepsilon)}_t\|_{1\to\infty}\leq a' \,t^{-d/2}
\label{etdy5.5}
\end{equation}
for all $t\geq1$.
These estimates are again uniform for $\varepsilon\in\langle0,1]$.

The estimates (\ref{etdy5.5}) convert to large time Gaussian bounds, with the distance associated 
with $C_\varepsilon$,  by Davies perturbation theory as in Proposition~\ref{padty4.1},
but with $\gamma=1$ and $\nu=0$.
Specifically one deduces that for all $\delta>0$ there exists an  $a'>0$ such that
\[
K^{(\varepsilon)}_t(x\,;y)
\leq a' \,  t^{-d/2} \,  
e^{-d_{C_\varepsilon}(x;y)^2 ((4+\delta)t)^{-1}}
\]
uniformly for all $t \geq 1$, $x,y \in \Ri^d$ and $\varepsilon\in \langle0,1]$. 
Finally, taking the limit $\varepsilon\to0$,   one obtains the upper bounds (\ref{etdy4.111}) on
$K^{(0)}$ by using Lemmas~\ref{leo4.1} and \ref{ladty4.2}.

The proof of the lower bounds is a repetition of the argument used to prove Corollary~\ref{ctdy3.3}.
Now one chooses $\psi$ to be the characteristic function of the ball $B_e(x_0\,;R\,t^{1/2})$ 
and uses the upper bounds
(\ref{etdy4.111}).
Moreover one uses the lower bound $d_C(x\,;y)\geq \|C\|^{-1/2}|x-y|$ to express the estimate of 
$|(\varphi, S^{(0)}_t\psi)|$
in terms of Euclidean parameters.
We omit the details.
\hfill$\Box$

\ruimte

Theorem~\ref{tadty411} implies that subellipticity and  local positivity gives
estimates 
\[
a'\,t^{-d/2}\leq \|K^{(0)}_t\|_\infty\leq a\,t^{-d/2}
\]
for all $t\geq1$.
Thus the asymptotic behaviour of the kernel is determined by the Euclidean dimension $d$
and is independent of the geometry related to the distance $d_C$.
This confirms the conclusions of Kusuoka and Stroock \cite{KuS}.
In fact Theorem~\ref{tadty411} can be applied directly to the class of operators
covered in Kusuoka and Stroock's main application Theorems~3.20 and 3.24.
It is possible to verify that their subellipticity assumption (3.21) implies our assumption
(\ref{eadty4.28}) and one can also verify that the kernels associated with their operators
are uniformly continuous.
Hence the local lower bounds assumed in  Theorem~\ref{tadty411} follow from the subellipticity by
Remarks~\ref{reo4.1} and \ref{reo4.2}.
These verifications will be contained in a separate article \cite{ER27}.

Although the bounds of Theorem~\ref{tadty411} verify the asymptotic behaviour suggested
by Kusuoka and Stroock they are weaker than the conclusions of these authors in two respects.
First the estimates do not give a Gaussian lower bound.
Secondly the estimates rely on an explicit assumption of  uniform  local positivity for small
$t$ of uniform continuity.
These features are, however, related and neither can be improved without
further assumption.
This will be  established by examples in Section~\ref{Sadty5}.
The problem arising in the case of degenerate operators is that the semigroup kernel is
not necessarily strictly positive.
In fact it  can take the value zero on sets  $(x,y)\in\Ri^d\times\Ri^d$ of non-zero measure.
This behaviour occurs if $\gamma\in\langle0,1/2]$.
It is possible that the kernel is strictly positive whenever $\gamma\in\langle 1/2,1\rangle$ and that the
for this range the   properties of the kernel resemble those found for  strongly elliptic 
operators , i.e., for operators with $\gamma=1$.
Alternatively  it could be relevant that on a large class manifolds for which the heat kernel 
satisfies Gaussian upper bounds the matching lower bounds are equivalent to H\"older
continuity of the kernel (see \cite{Cou4} and the extensive list of references therein).

\section{Subelliptic resolvent estimates}\label{Sadty4.1}

In this section we use the off-diagonal bounds of Section~\ref{Sadty3} to establish pointwise
on-diagonal lower bounds for the kernel of a high power of the resolvent.
The bounds are  more efficient than the earlier bounds since they are position dependent.
Moreover they can be inverted to give lower bounds on the Euclidean volume of the balls 
$B_C$ defined by the
quasidistance~$d_C$.

First we derive a  statement for strongly elliptic operators which will be applied to the 
approximants $H_\varepsilon$.
We use the notation $K_S$ for the distribution kernel of a bounded operator $S$.

\begin{thm}\label{tadty511}
For all $m \in \Ni$ with $4m > d$ there exists an $a > 0$ such that for any
strongly elliptic operator $H$ with measurable coefficients 
\[
a \, |B_C(x\,;r)| \geq K_{(I + r^2 H)^{-2m}}(x\,;x)^{-1}
\geq \|(I + r^2 H)^{-m}\|_{2\to\infty}^{-2}
\]
for all $x \in \Ri^d$ and $r > 0$.
\end{thm}

Note that $K_{(I + r^2 H)^{-2m}}$ is continuous by standard estimates for
strongly elliptic operators if $4m > d$.
Therefore its on-diagonal value  is well-defined.

The proof of Theorem~\ref{tadty511} requires an extension of Lemma~\ref{lod10}.

\begin{lemma}\label{ladty501}
If $H$ is strongly elliptic and $\Psi \in {\cal S}(\Ri)$ is an  even function with 
$\supp \Psi\subseteq [-1,1]$ then the kernel $K_{\widehat\Psi(r H^{1/2})}$ is continuous
for all $r > 0$ and 
\[
K_{\widehat\Psi(r H^{1/2})}(x\,;y) = 0
\]
for all $x,y \in \Ri^d$ with $r < d_C(x\,;y)$.
\end{lemma}
\proof\
For all $m \in \Ni$ define $\Phi_m \in {\cal S}(\Ri)$ by
$\Phi_m(\lambda) = \widehat\Psi(\lambda) (1 + \lambda^2)^m$.
We first show that the distributional kernel 
$K_{\widehat \Psi(r H^{1/2})}$ of $\widehat \Psi(r H^{1/2})$ is continuous.
Note that 
\[
\widehat \Psi(r H^{1/2})
= (I + r^2 H)^{-m} \, \Phi_{2m}(r H^{1/2}) \, (I + r^2 H)^{-m}
\;\;\; .  \]
Then $\Phi_{2m}(r H^{1/2})$ maps $L_2$ into $L_2$.
Moreover, if $4m > d$ then the resolvent
$(I + r^2 H)^{-m}$ maps $L_2$ into $L_\infty$ and $L_1$ into $L_2$.
So $\widehat \Psi(r H^{1/2})$ maps $L_1$ into $L_\infty$ and the 
distributional kernel $K_{\widehat \Psi(r H^{1/2})}$ is a bounded function.
Since $\Phi_{2m} \in {\cal S}(\Ri)$ the 
same applies to the kernel $K_{ \Phi_{2m}(r H^{1/2})}$
of $ \Phi_{2m}(r H^{1/2})$.
But if $2m > d$ then the kernel $K_{(I + r^2 H)^{-m}}$ of $(I + r^2 H)^{-m}$ 
is continuous, 
$\sup_{x \in \Ri^d} \int_{\Ri^d} dy \, |K_{(I + r^2 H)^{-m}}(x\,;y)| < \infty$
and there are $c,\nu > 0$ such that 
\[ 
\int_{\Ri^d} dy \, |K_{(I + r^2 H)^{-m}}(x_1\,;y) - K_{(I + r^2 H)^{-m}}(x_2\,;y)|
\leq c \, |x_1 - x_2|^\nu
\]
for all $x_1,x_2 \in \Ri^d$ as a result of standard estimates for
strongly elliptic operators.
Since
\[
K_{\widehat \Psi(r H^{1/2})}(x\,;y)
= \int dz_1 \int dz_2 \, K_{(I + r^2 H)^{-m}}(x\,;z_1) \, 
K_{ \Phi_{2m}(r H^{1/2})}(z_1\,;z_2) \, 
K_{(I + r^2 H)^{-m}}(z_2\,;y)
\]
for all $x,y \in \Ri^d$ it follows that $K_{\widehat \Psi(r H^{1/2})}$ is H\"older continuous
and in particular continuous.

Let $x_1,x_2 \in \Ri^d$ and suppose that $d_C(x_1\,;x_2) > r$.
Set $\varepsilon = 2^{-1} (d_C(x_1\,;x_2) - r)$.
Then Lemma~\ref{lod10} states that 
\[
\int dx \int dy \, \overline{\varphi_1(x)} \, 
K_{\widehat \Psi(r H^{1/2})}(x\,;y) \, \varphi_2(y)
= (\varphi_1, \widehat \Psi(r H^{1/2}) \varphi_2)
= 0
\]
for all $\varphi_1 \in L_2(B_C(x_1\,;\varepsilon))$ and 
$\varphi_2 \in L_2(B_C(x_2\,;\varepsilon))$.
Since $H$ is strongly elliptic the metric~$d_C$ is equivalent to the 
Euclidean metric $d_e$ on $\Ri^d$.
Hence there is a $\delta > 0$ such that $B_e(x_1\,;\delta) \subset B_C(x_1\,;\varepsilon)$
and $B_e(x_2\,;\delta) \subset B_C(x_2\,;\varepsilon)$, where $B_e$ denotes the 
Euclidean ball.
Then 
\[
\int dx \int dy \, \overline{\varphi_1(x)} \, 
K_{\widehat \Psi(r H^{1/2})}(x\,;y) \, \varphi_2(y)
= 0
\]
for all $\varphi_1 \in L_2(B_e(x_1\,;\delta))$ and 
$\varphi_2 \in L_2(B_e(x_2\,;\delta))$.
Since $K_{\widehat \Psi(r H^{1/2})}$ is continuous this implies that 
$K_{\widehat \Psi(r H^{1/2})}(y_1\,;y_2) = 0$ for all 
$y_1 \in B_e(x_1\,;\delta)$ and $y_2 \in B_e(x_2\,;\delta)$.
In particular $K_{\widehat \Psi(r H^{1/2})}(x_1\,;x_2) = 0$.\hfill$\Box$

\ruimte

\noindent
{\bf Proof of Theorem~\ref{tadty511}} \hspace{5pt}\
Fix $\Psi \in {\cal S}(\Ri)$ even  with 
$\supp \Psi\subseteq [-1,1]$ and $\int \Psi = 1$.

It follows from Corollary~\ref{cadt213}, applied with $\Phi=\widehat \Psi$, that 
$\widehat \Psi(r H^{1/2}) \, \one = \widehat \Psi(0) \, \one = \one$
for all $r > 0$.
Let $x \in \Ri^d$ and $r > 0$.
Then the support property of Lemma~\ref{ladty501} and the Cauchy--Schwarz inequality imply  that 
\begin{equation}
1 
= \int dy \, K_{\widehat \Psi(r H^{1/2})}(x\,;y) 
 \leq  |B_C(x\,;r)|^{1/2} 
\Big( \int dy \, |K_{\widehat \Psi(r H^{1/2})}(x\,;y)|^2 \Big)^{1/2} 
\label{epadt301;1}
\;\;\; .  
\end{equation}
Let $m \in \Ni$ with $4m > d$.
Define $\Phi_m \in {\cal S}(\Ri)$ by
$\Phi_m(\lambda) = \widehat\Psi(\lambda) (1 + \lambda^2)^m$ as before.
Then 
\[
\widehat \Psi(r H^{1/2})
= (I + r^2 H)^{-m} \, \Phi_m(r H^{1/2})
\]
and $\|\Phi_m(r H^{1/2})\|_{2\to2}=\|\Phi_m\|_\infty$.
Now if $S$ is bounded from $L_2(\Ri^d)$ to $L_\infty(\Ri^d)$ and $T$ is bounded from
$L_2(\Ri^d)$ to $L_2(\Ri^d)$ then
\[
\int dy\,|K_{ST}(x\,;y)|^2\leq \|T\|^2_{2\to2}\int dy\,|K_S(x\,;y)|^2
\;\;\;.
\]
Therefore applying this estimate to (\ref{epadt301;1}) with $S=(I+r^2H)^{-m}$ and $T=\Phi_m(r H^{1/2})$
one finds
\begin{eqnarray*}
1 & \leq & |B_C(x\,;r)| \,\|\Phi_m\|_\infty^2  \int dy \, |K_{(I+r^2H)^{-m}}(x\,;y)|^2 \\[5pt]
&=& |B_C(x\,;r)| \,\|\Phi_m\|_\infty^2 \,K_{(I+r^2H)^{-2m}}(x\,;x)
\leq  |B_C(x\,;r)| \,\|\Phi_m\|_\infty^2 \,\|(I + r^2 H)^{-m}\|_{2\to\infty}^2
\end{eqnarray*}
where the second  relation uses $K_S*K_T=K_{ST}$.
\hfill$\Box$

\ruimte

Note that the statement of the proposition could be inverted to give
\[
K_{(I+tH)^{-2m}}(x\,;x)\geq a_m^{-1}\, |B_C(x\,;t^{1/2})|^{-1}
\;\;\;,
\]
i.e., one has an on-diagonal lower bound for the kernel with the anticipated spatial dependence.

Next we consider lower bounds 
on the volume of the balls associated with a subelliptic operator.
In the subelliptic situation a new phemonenon of separation occurs.
This will be  discussed in detail
in Section~\ref{Sadty5}.
It is possible to have subspaces $L_2(\Omega)$ of $L_2(\Ri^d)$ which are invariant under
$S^{(0)}_t$ for all $t>0$.
The following result is adapted to this situation.

\begin{thm} \label{tadt302}
Assume the viscosity operator $H_0$ is subelliptic of order $\gamma$.
Further assume that there is a non-empty open subset $\Omega\subseteq \Ri^d$
such that $S^{(0)}_tL_2(\Omega)\subseteq L_2(\Omega)$.

Then there are $a,R > 0$ such that
\begin{equation}
|B_C(x\,;r)\cap\Omega| \geq a \, r^{d/\gamma}
\label{eadt302;2}
\end{equation} 
for all $x \in \Omega$ and $r \in \langle0,R]$.
The values of $a$ and $R$ depend only on the subellipticity parameters $\mu$, $\nu$ and $\gamma$
and are independent of $\Omega$.
\end{thm}
\proof\
Let $x \in \Omega$ and $r > 0$.
For all $\varepsilon > 0$ let 
$H_\varepsilon$ denote the strongly elliptic approximants to $H_0$
and let $\Psi \in {\cal S}(\Ri)$ be even  with 
$\supp \Psi\subseteq [-1,1]$ and $\int \Psi = 1$.
Since ${\widehat \Psi}(0)=1$ and ${\widehat \Psi}(rH_\varepsilon^{1/2})$ is symmetric
it follows from Corollary~\ref{cadt213}, applied with $\Phi=\widehat \Psi$, that 
\[
(\one,\varphi)=(\one,{\widehat \Psi}(r H_\varepsilon^{1/2}) \varphi)
\]
for all $\varphi\in C_c^\infty(\Ri^d)$.

Next since $\Omega$ is open there exists an $s_0 \in \langle0,1]$ such that 
$B_{C_1}(x\,;s_0)\subset\Omega$.
Let $\delta \in \langle0,1]$ and $s \in \langle0,s_0]$.
Then $B_{C_\delta}(x\,;s) \subseteq B_{C_1}(x\,;s_0) \subset\Omega$.
Let $\varepsilon\in\langle0,\delta]$.
If $x_0,y_0 \in \Ri^d$ and $d_{C_\delta}(x_0\,;y_0)>r$ then 
$d_{C_\varepsilon}(x_0\,;y_0) \geq d_{C_\delta}(x_0\,;y_0)>r$
and $K_{{\widehat\Psi}(rH_\varepsilon^{1/2})}(x_0\,;y_0)=0$ by Lemma~\ref{ladty501}.
But there exists a positive, non-zero $\varphi\in C_c^\infty(\Ri^d)$ with 
$\supp\varphi\subset B_{C_\delta}(x\,;s)$.
Then $\supp{\widehat \Psi}(r H_\varepsilon^{1/2}) \varphi\subseteq B$
where for brevity we have set $B= B_{C_\delta}(x\,;r+s)$.
Therefore
\[
(\one,\varphi)=(\one_B,{\widehat \Psi}(r H_\varepsilon^{1/2}) \varphi)
\;\;\; .  \]
Since the ball $B$ is relatively compact the characteristic function $\one_B$
is an $L_2$-function and 
since $H_\varepsilon$ converges in the strong resolvent sense to $H_0$ it follows that 
${\widehat\Psi}(rH_\varepsilon^{1/2})$ converges strongly to ${\widehat\Psi}(rH_0^{1/2})$
on $L_2(\Ri^d)$ by \cite{RS1} Theorem~VIII.20.
Therefore one deduces that 
\[
(\one,\varphi)=(\one_B,{\widehat \Psi}(r H_0^{1/2}) \varphi)
\;\;\;.
\]
Since  $\varphi \in L_2(\Omega)$ and $S^{(0)}_t$
leaves $L_2(\Omega)$ invariant it follows that 
${\widehat \Psi}(r H_0^{1/2}) \varphi \in L_2(\Omega)$.
Consequently
\begin{eqnarray*}
\|\varphi\|_1=(\one,\varphi)
=(\one_{B\cap\Omega},{\widehat \Psi}(r H_0^{1/2}) \varphi)
\leq\|\one_{B\cap\Omega}\|_2\,\|{\widehat \Psi}(r H_0^{1/2})\|_{1\to2}\,\|\varphi\|_1
\end{eqnarray*}
from which one concludes that
\[
1\leq |B_{C_\delta}(x\,;r+s)\cap\Omega|^{1/2}\,\|{\widehat \Psi}(r H_0^{1/2})\|_{1\to2}
\]
for all $\delta\in\langle0,1]$ and $s \in \langle0,s_0]$.
But $|B_{C_\delta}(x\,;r+s)\cap\Omega|$ decreases as $\delta$ and $s$ decrease to zero
and
\[
B_{C_\delta}(x\,;r+s)\subset B_{C_1}(x\,;r+1)
\]
for all $\delta\in\langle0,1]$ and $s \in \langle0,s_0]$.
Moreover $|B_{C_1}(x\,;r+1)|<\infty$.
Therefore
\begin{equation}
1\leq |B_C(x\,;r)\cap\Omega|^{1/2}\,\|{\widehat \Psi}(r H_0^{1/2})\|_{1\to2}
\label{etadt302;1}
\end{equation}
by Lemma~\ref{ladty4.2}.

Next let $m \in \Ni$.
If $\Phi_m(\lambda) = \widehat\Psi(\lambda) (1 + \lambda^2)^m$ as before then one has
\[
\|{\widehat \Psi}(rH_\varepsilon^{1/2})\varphi\|_2\leq \|\Phi_m\|_\infty\|(I+r^2H_\varepsilon)^{-m}\varphi\|_2
\]
for all $\varphi\in L_2(\Omega)$, $m \in \Ni$ and $\varepsilon > 0$ 
 as in the proof of Theorem~\ref{tadty511}.
Then in the limit $\varepsilon\to0$ one deduces that 
\begin{equation}
\|{\widehat \Psi}(r H_0^{1/2})\|_{1\to2}\leq \|\Phi_m\|_\infty\|(I+r^2H_0)^{-m}\|_{1\to2}
\label{etadt302;3}
\end{equation}
for all $m \in \Ni$.

Finally it follows from the proof of Proposition~\ref{padty4.1} that 
there is an $a > 0$, depending only on $\gamma$ and $d$, such that 
\[
\|S^{(0)}_t\|_{1\to2}
\leq a \,(\mu t)^{-d/(4\gamma)}e^{\nu t}
\]
for all $ t > 0$.
Then one estimates
\begin{eqnarray*}
\|(I+r^2H_0)^{-m}\|_{1\to2}
&\leq& ((m-1)!)^{-1}\int^\infty_0dt\,e^{-t}t^{m-1}
\|S^{(0)}_{r^2t}\|_{1\to2}\\[5pt]
&\leq& a\,((m-1)!)^{-1}\int^\infty_0dt\,e^{-t}t^{m-1}
(\mu r^2t)^{-d/(4\gamma)}e^{\nu r^2t}
\end{eqnarray*}
and the integral is finite if $r^2 < \nu^{-1}$ and $m>d/(4\gamma)$.
Now fix $m \in \Ni$ with $m>d/(4\gamma)$ and set $R = (2 \nu)^{-1/2}$.
Then there is an $a' > 0$ such that 
\begin{equation}
\|(I+r^2H_0)^{-m}\|_{1\to2}
\leq a' \,r^{-d/(2\gamma)}
\label{etadt302;2}
\end{equation}
uniformly for all $r \in \langle0,R]$.
Therefore the theorem follows by a combination of (\ref{etadt302;1}), (\ref{etadt302;3}) and 
(\ref{etadt302;2}).\hfill$\Box$

\ruimte

The statement of the theorem is related to Theorem~1 of Fefferman and Phong \cite{FP}.
The latter result establishes for operators with smooth coefficients  that subellipticity
gives a local comparison, $d_C(x\,;y)\leq a\,|x-y|^\gamma$ for $x,y\in\Ri^d$ with $|x|,|y|\leq1$, of the distance $d_C$ and
the subelliptic distance and the Euclidean distance. 
The statement of Theorem~\ref{tadt302} is of a similar nature.
It is weaker insofar it only compares the volume of balls but it is stronger insofar it is global and 
valid for operators with measurable coefficients.
In fact we next show that it gives the Fefferman--Phong result in one-dimension without 
any smoothness requirements.

In one-dimension  $H$ is  formally given by $H=-d\,c\,d$, with $d=d/dx$ 
and $0 \leq c \in L_\infty(\Ri)$, and 
the strongly elliptic approximants  by $H_\varepsilon=-d\,(c+\varepsilon)\,d$.
The distance $d_{C_\varepsilon}$ is easily computed to be
\begin{equation}
d_{C_\varepsilon}(x\,;y)=\int^y_xdz\,(c(z)+\varepsilon)^{-1/2}
\label{eadt401}
\end{equation}
for all $x,y\in \Ri$ with $y>x$ and 
$d_C(x\,;y)=\lim_{\varepsilon\to0}d_{C_\varepsilon}(x\,;y)$ by 
Lemma~\ref{ladty4.2}.
Now the volume estimates of Theorem~\ref{tadt302} allow one to deduce that $d_C$ is 
finite-valued.

\begin{thm}\label{tadty512}
Assume $d=1$ and $H_0$ is subelliptic of order $\gamma$.
Then $c(x) > 0$ for almost all $x \in \Ri$ and 
\begin{equation}
d_C(x\,;y) = \int^y_x dz\,c(z)^{-1/2}
\label{eadt412}
\end{equation}
for all $x,y \in \Ri$ with $x < y$.
Moreover, there exist $a_1,a_2>0$ such that
\begin{equation}
a_1\,|x-y|\leq d_C(x\,;y)\leq a_2\,( |x-y|^\gamma \vee |x-y| )
\label{eadt41;1}
\end{equation}
for all $x,y \in \Ri$.
In particular $d_C$ is finite-valued.
\end{thm}
\proof\
The proof is based upon the volume estimates of Theorem~\ref{tadt302} applied with $\Omega = \Ri$.
These are valid since we assume $H_0$ to be subelliptic.

First one has $d_C(x\,;y)\geq \|c\|_\infty^{-1/2}\,|x-y|$.
Secondly, since $B_{C_\varepsilon}(x\,;r)\supseteq B_C(x\,;r)$ the volume bounds give
\[
|B_{C_\varepsilon}(x\,;r)|\geq a\, r^{1/\gamma}
\]
for all  $r \in \langle0,R]$ uniformly for $\varepsilon>0$.
Now $d_{C_\varepsilon}(x\,;y)<\infty$ for all $x,y\in\Ri$.
But, assuming always that $x < y$,
there is a $z \in \langle x,y\rangle$ such that 
$d_{C_\varepsilon}(x\,;z) = d_{C_\varepsilon}(z\,;y)$.
Then $B_{C_\varepsilon}(z\,;r) = \langle x,y\rangle$ where 
$r = 2^{-1} d_{C_\varepsilon}(x\,;y)$.
So if $d_{C_\varepsilon}(x\,;y) \leq 2R$ then 
\[
|y-x| 
= |B_{C_\varepsilon}(z\,;r)|
\geq a \, r^{1/\gamma}
= a\, (2^{-1} d_{C_\varepsilon}(x\,;y))^{1/\gamma}
\]
and 
\begin{equation}
d_{C_\varepsilon}(x\,;y) \leq a' \, |x-y|^\gamma
\label{ewip1;12}
\end{equation}
where $a' = 2 a^{-\gamma}$ is independent of $\varepsilon$.

Now let $x,y \in \Ri$ and suppose $2 a' |x-y|^\gamma \leq R$.
We claim that 
\begin{equation}
d_{C_\varepsilon}(x\,;y)
\leq 2a' \, |x-y|^\gamma
\label{ewip1;13}
\end{equation}
for all $\varepsilon > 0$.
Indeed, if $\varepsilon \in \langle1,\infty\rangle$ is large then
$d_{C_\varepsilon}(x\,;y) \leq \varepsilon^{-1/2} |x-y| \leq R$.
Hence $d_{C_\varepsilon}(x\,;y) \leq a' \, |x-y|^\gamma$ by (\ref{ewip1;12}).
If (\ref{ewip1;13}) is not valid for all $\varepsilon > 0$ then 
set 
\[
\varepsilon_0
= \sup \{ \varepsilon > 0 : d_{C_\varepsilon}(x\,;y) > 2a' \, |x-y|^\gamma \}
\;\;\; .  \]
By (\ref{eadt401}) the function $\varepsilon \mapsto d_{C_\varepsilon}(x\,;y)$ is continuous and 
decreasing.
Therefore $d_{C_{\varepsilon_0}}(x\,;y) = 2a' \, |x-y|^\gamma \leq R$
and $d_{C_{\varepsilon_0}}(x\,;y) \leq a' \, |x-y|^\gamma$
by (\ref{ewip1;12}).
This is a contradiction.

So $d_{C_\varepsilon}(x\,;y) \leq 2a' \, |x-y|^\gamma$ for all $x,y \in \Ri$ 
and $\varepsilon > 0$ if $2 a' |x-y|^\gamma \leq R$.
Then 
\[
d_C(x\,;y) = \sup_{\varepsilon > 0} d_{C_\varepsilon}(x\,;y) \leq 2a' \, |x-y|^\gamma
\]
for all $x,y \in \Ri$ with $2 a'\,|x-y|^\gamma \leq R$ by 
Lemma~\ref{ladty4.2}.
Hence by the triangle inequality there is an $a'' > 0$ such that 
$d_C(x\,;y) \leq a'' \, |x-y|^\gamma$ for all $x,y \in \Ri$ with $|x-y| \leq 1$.
Therefore, if $|x-y| \geq 1$ it follows again from the triangle inequality that 
$d_C(x\,;y) \leq 2 a'' \, |x-y|$.
Combining these bounds one deduces that 
\[
d_C(x\,;y)
\leq 2 a''( |x-y|^\gamma \vee |x-y| )
\]
for all $x,y \in \Ri$.
This proves the upper bounds of (\ref{eadt41;1}) with $a_2 = 2 a''$.

Next  since $d_{C_\varepsilon}(x\,;y)$ is given by  (\ref{eadt401})
and 
$\lim_{\varepsilon\to0}d_{C_\varepsilon}(x\,;y) = d_C(x\,;y)$ by Lemma~\ref{ladty4.2}.
But $d_C$ 
is finite valued
and so  it follows that 
$c(x) > 0$ for almost every $x \in \Ri$.
But then
\[
d_C(x\,;y)
= \lim_{\varepsilon\to0}d_{C_\varepsilon}(x\,;y)
=\int^y_x dz\,c(z)^{-1/2}
\]
for all $x,y \in \Ri$ with $x < y$.\hfill$\Box$

\ruimte

A much stronger conclusion is valid if the coefficient $c$ is twice-differentiable.
\begin{prop}\label{padty513}
Let $d=1$. Assume $c\geq0$,  $c \in C_b^2(\Ri)$ and the corresponding viscosity operator
$H_0$ is subelliptic.
Then $H_0$  is strongly elliptic.
 \end{prop}
\proof\
It follows from the above argument  $c > 0$ almost everywhere and one has
\[
c(y) - 2 c(x)
\leq c(y) + c(2x-y) - 2 c(x)
= \int_0^{y-x} dt \, \int_{x-t}^{x+t} ds \, c''(s)
\leq \|c''\|_\infty (y-x)^2
\]
for all $x,y \in \Ri$.
Now let $x \in \Ri$. 
Then for all $y \in \langle x,x+1]$ one has by Cauchy--Schwarz

\begin{eqnarray*}
y-x
& = & \int_x^y c^{-1/4} \, c^{1/4}  \\[5pt]
& \leq & \Big( \int_x^y c^{-1/2} \Big)^{1/2}
         \Big( \int_x^y  c^{1/2} \Big)^{1/2}  \\[5pt]
& \leq & d_C(x\,;y)^{1/2}
  \Bigg( \Big( \int_x^y  \one \Big)^{1/2}
         \Big( \int_x^y c \Big)^{1/2} \Bigg)^{1/2}  \\[5pt]
& \leq & a_2^{1/2} |y-x|^{\gamma/2} \, |y-x|^{1/4} 
            \Big( \int_x^y c \Big)^{1/4}
\end{eqnarray*}
where $a_2 > 0$ is as in (\ref{eadt41;1}).
Therefore
\begin{eqnarray*}
a_2^{-2} |y-x|^{3 - 2\gamma}
& \leq & \int_x^y dt \, c(t)  \\[5pt]
& \leq & \int_x^y dt \Big( 2 c(x) + \|c''\|_\infty \, (t-x)^2 \Big) \\[5pt]
& = & 2 c(x) \, |y-x| + 3^{-1} \|c''\|_\infty \, |y-x|^3
\;\;\; .  
\end{eqnarray*}
Rearranging gives
\[
c(x)
\geq 2^{-1} |y-x|^{2-2\gamma} 
       \Big( a_2^{-2} - 3^{-1} \|c''\|_\infty \, |y-x|^{2\gamma} \Big)
\;\;\; .  \]
Now choose $y = x + (3 /(2 a_2^2 \|c''\|_\infty + 3))^{1/(2\gamma)}$.
Then
\[
c(x) \geq (4 a_2^2)^{-1} (3/ (2 a_2^2 \|c''\|_\infty + 3))^{(1-\gamma)/\gamma}
\]
and $H_0$ is strongly elliptic.\hfill$\Box$

\ruimte

In the context of operators with smooth coefficients Fefferman and Phong derived a converse 
 statement that local comparability of $d_C$ and the Euclidean distance, i.e., estimates
of the form (\ref{eadt41;1}), imply subellipticity.
But no such general statement is possible for operators with measurable coefficients.
We will give counterexamples in Section~\ref{Sadty5}.
The examples even have coefficients in $C^{2\gamma}_b$ for $\gamma\in \langle0,1\rangle$.
But we conclude this section with some simple examples of subelliptic operators.

\begin{exam}\label{exadty4.1}
Let   $(c_{ij}) \geq  (c \,\delta_{ij})$ with 
\[
c(x)= \Big( \frac{|x|^2}{1+|x|^2} \Big)^\delta
\]
and $\delta\in[0,1\rangle$.
If $ \delta<d/2$ then we  prove that the corresponding $H_0$ is subelliptic of order $1-\delta$.

First suppose $d \geq 3$.
Then one has  the elementary quadratic form inequality $\Delta\geq\sigma\,|x|^{-2}$ 
with $\sigma=(d-2)^2/4$
(see, for example, \cite{Kat1} Remark~VI.4.9a and (VI.4.24), or \cite{RS2} Lemma on page~169).
The inequality immediately implies  that $c \geq a_0\,(I+\Delta)^{-\delta}$ as quadratic forms
for a suitable $a_0>0$.
Then it follows from Proposition~\ref{padty215} that 
\begin{equation}
H_0\geq a_0\,\Delta(I+\Delta)^{-\delta}=a_0\, \Delta^{1-\delta}(\Delta (I+\Delta)^{-1} )^\delta
\;\;\;.
\label{eadt303}
\end{equation}
But since $\delta\in[0,1\rangle$ one has
\[
H_0\geq a_0\,(\varepsilon(1+\varepsilon)^{-1})^\delta
\int^\infty_\varepsilon dE_\Delta(\lambda)\,\lambda^{1-\delta}
\geq  a_0\,(\varepsilon (1+\varepsilon)^{-1} )^\delta\,\Big(\Delta^{1-\delta}
-\varepsilon^{1-\delta}\,I\Big)
\]
for all $\varepsilon>0$ where $E_\Delta$ denotes the spectral family of $\Delta$.
Thus the operator $H_0$ is subelliptic of order $1-\delta$.

A similar conclusion holds  for $d=2$ with $\delta\in[0,1\rangle$ and $d=1$ 
with $\delta\in[0,1/2\rangle$
by the following fractional version of the foregoing argument.

It follows by a general result of Strichartz on multipliers on Sobolev spaces 
\cite{Stri}, Theorem 3.6,
that $|x|^{-\delta}\Delta^{-\delta/2}$ is bounded on $L_2(\Ri^d)$ if $\delta\in[0,d/2\rangle$
 (see \cite{RS2}, Chapter~IX, Exercise~39(b)).
Therefore 
$\Delta^{-\delta/2}|x|^{-2\delta}\Delta^{-\delta/2}$ is bounded and this means that 
there is a $\sigma>0$ such that $\Delta^\delta\geq \sigma\,|x|^{-2\delta}$
in the sense of quadratic forms.
This is a fractional version of the foregoing estimate and it again gives a  bound 
$c \geq a_0\,(I+\Delta)^{-\delta}$ as forms
for a suitable $a_0>0$.
Hence    the estimate 
(\ref{eadt303}) is now valid under the restriction
  $\delta\in[0,d/2\rangle$
and the subellipticity estimate (\ref{eadty4.28}) again follows for $\delta<1\wedge (d/2)$ 
by the foregoing spectral argument.
Thus  if $d=2$  the estimate is valid for all $\delta\in[0,1\rangle$
but if $d=1$ it is only established for $\delta\in[0,1/2\rangle$.   
The situation  for $d=1$ and $\delta\in[1/2,1\rangle$
is more complicated.
It will be discussed in detail in Examples~\ref{exadty5.1} and \ref{xadty606}.
\hfill$\Box$
\end{exam}

  Although the coefficients in these examples are only degenerate at the single point $x=0$
it is easy to construct examples with a finite number of degeneracies with different orders 
of degeneracy.

\begin{exam}\label{exadty4.2}
Let  $x_1,\ldots,x_n$ be distinct points in $\Ri^d$ and 
$\delta_1,\ldots,\delta_n\in\langle0,1\rangle$.
Set $\underline \delta=\min \delta_i$,
$\overline \delta=\max \delta_i$ and define $d \colon \Ri^d \to [0,\infty\rangle$ by
$d(x)=\min_{1\leq i\leq n}|x-x_i|^{\delta_i/{\overline \delta}}$.
Now if $d\geq 3$  it follows from the bounds $\Delta\geq\sigma\,|x|^{-2}$  and  
translation invariance that
 $\Delta \geq\sigma\,|x-x_i|^{-2}$.
Therefore one finds straightforwardly that there is a $\sigma' > 0$ such that 
\[
d(x)^{-2} 
= \max_{1\leq i\leq n} |x-x_i|^{-2\delta_i/{\overline \delta}}
\leq \sum_{i=1}^n |x-x_i|^{-2\delta_i/{\overline \delta}}
\leq \sum_{i=1}^n \sigma^{-\delta_i/{\overline \delta}} \Delta^{{\delta_i/{\overline \delta}}}
\leq \sigma' (\Delta+\Delta^{{\underline \delta}/{\overline \delta}})
\;\;\; .  \]
Now consider  operators with  
 $(c_{ij}) \geq c\,\delta_{ij}$ where 
\[
c(x)= \Big( \frac{d(x)^2}{1+d(x)^2} \Big)^\delta
\]
and $\delta \in [0,1\rangle$.
It follows as above that there are $a,a'>0$ such that 
\[
H_0\geq a\,\Delta(I+\Delta+\Delta^{{\underline \delta}/{\overline \delta}})^{-\delta}
\geq a'\, \Delta(I+\Delta)^{-\delta}
\;\;\;.
\]
Then by spectral theory there are $\mu,\nu>0$ such that 
\[
H_0\geq\mu\,\Delta^{1-\delta}-\nu\,I
\]
and $H_0$ is subelliptic of order $1-\delta$.
Similarly if $d\leq 2$  the fractional bounds 
$\Delta^\delta\geq \sigma\,|x|^{-2\delta}$ imply that 
$(\Delta+\Delta^{{\underline \delta }/{\overline \delta} })^\delta\geq \sigma'\,d(x)^{-2\delta}$.
Hence  one can establish that $H_0$ satisfies the subelliptic condition 
whenever $\delta<1\wedge d/2$.
In the next section (see Examples~\ref{exadty5.1} and \ref{xadty606}) we return to the discussion of the situation  for  $d=1$ and 
$ \delta\in[ 1/2,1 \rangle$.
\hfill$\Box$
\end{exam}

\section{Separation properties}\label{Sadty5}

The foregoing properties of elliptic and subelliptic operators are direct analogues
of similar properties of strongly elliptic operators.
The principal difference is the replacement of the Euclidean distance by the distance
$d_C$.
But now we examine a phenomenon which has no analogue for strongly elliptic operators,
the phenomenon of separation either partial or complete.
Degeneracy of the coefficients  can lead to the system factoring into independent subsystems,
i.e, there is a complete separation.
It is also possible to have an incomplete separation but we will not examine
this behaviour.
These phenomena do not require any particular pathological property of the coefficients 
and can occur even if the coefficients are nearly $C^2$ and the operator is subelliptic.

In the sequel we shall need the following simple lemma.

\begin{lemma} \label{ladty601}
Let $\Omega$ be a measurable subset of $\Ri^d$.
Let $S$ be a bounded self-adjoint operator on $L_2(\Ri^d)$ which  extends
to a bounded operator on $L_p(\Ri^d)$ for all $p \in [1,\infty]$.
Suppose $S$ is positivity preserving.
Then the following are equivalent.
\begin{tabel}
\item \label{ladty601-1}
There exists a $p \in [1,\infty]$ such that $S L_p(\Omega) \subseteq L_p(\Omega)$.
\item \label{ladty601-4}
There exists a $c > 0$ such that $S \one_\Omega \leq c\,  \one_\Omega$.
\end{tabel}
If $S \one = \one$ then $S \one_\Omega = \one_\Omega$ in Statement~{\rm \ref{ladty601-4}}.
\end{lemma}
\proof\
If Statement~\ref{ladty601-1} holds for some $p \in [1,\infty]$ then by a density 
argument   it is valid for  $p=2$. 
Next $L_2(\Omega)$ is a closed subspace of
 $L_2(\Ri^d)$ and $S$ is self-adjoint.
Hence by a standard argument $L_2(\Omega)$
is invariant for $S$   if and only if $S$ commutes with the orthogonal projection
on~$L_2(\Omega)$
\begin{equation}
SM_\Omega=M_\Omega S\;\;\;,
\label{eadty6.138}
\end{equation}
 where
$M_\Omega(\varphi)=\one_\Omega \varphi$.
Then, by another density argument, (\ref{eadty6.138}) holds also
on  $L_\infty$ and 
$S \one_\Omega = SM_\Omega \one = M_\Omega S \one =
\one_\Omega S \one \leq \|S \one \|_\infty \one_\Omega$. 
Thus
Statement~\ref{ladty601-4} holds 
with $c = \|S \one\|_\infty$.
In particular, if $S \one = \one$ then $S \one_\Omega = \one_\Omega S \one = \one_\Omega $.
Finally, if Statement~\ref{ladty601-4} is valid then for all real $\varphi \in L_2(\Omega)$ one has 
\[
|S \varphi|
\leq S |\varphi|
\leq S ( \|\varphi\|_\infty \one_\Omega)
\leq c \,\|\varphi\|_\infty \one_\Omega
\]
since $S$ is  positivity preserving.
Therefore $S \varphi \in L_\infty(\Omega)$.
Then it easily follows that Statement~\ref{ladty601-1} is valid for $p=\infty$.
This completes the proof of  the lemma.\hfill$\Box$

\begin{cor} \label{cadty602}
Let $S$ be a bounded self-adjoint operator on $L_2(\Ri^d)$ which extends
to a bounded operator on $L_p(\Ri^d)$ for all $p \in [1,\infty]$.
Suppose $S$ is positivity preserving.
Then 
\[
\{ \Omega \subseteq \Ri^d : \Omega \mbox{ is measurable and }
       S L_2(\Omega) \subseteq L_2(\Omega) \}
\]
is a $\sigma$-algebra.
\end{cor}

Separation of a semigroup of operators associated to a Dirichlet form can be characterized in terms of the 
Dirichlet form.

\begin{lemma} \label{ladty616}
Let $k$ be a Dirichlet form and $T$ the associated semigroup.
Let $\Omega$ be a measurable subset of $\Ri^d$.
The following are equivalent.
\begin{tabel} 
\item \label{ladty616-1}
$T_t L_2(\Omega) \subseteq L_2(\Omega)$ for all $t > 0$.
\item \label{ladty616-2}
For all $\varphi \in D(k)$ one has $\varphi \one_\Omega \in D(k)$ and 
\begin{equation}
k(\varphi) = k(\varphi \one_\Omega) + k(\varphi \one_{\Omega^{\rm c}})
\;\;\; .  
\label{eladty616;1}
\end{equation}
\end{tabel}
\end{lemma}
\proof\
\ref{ladty616-1}$\Rightarrow$\ref{ladty616-2}.
Let $\varphi \in D(k)$.
Then 
\begin{equation}
t^{-1} (\varphi, (I - T_t) \varphi)
= t^{-1} (\varphi \one_\Omega, (I - T_t) \varphi \one_\Omega)
   + t^{-1} (\varphi \one_{\Omega^{\rm c}}, (I - T_t) \varphi \one_{\Omega^{\rm c}})
\label{eladty616;2}
\end{equation}
for all $t > 0$. Moreover, both terms on the right hand side of (\ref{eladty616;2}) are 
positive.
Hence 
\[
\sup_{t > 0} t^{-1} (\varphi \one_\Omega, (I - T_t) \varphi \one_\Omega)
\leq \sup_{t > 0} t^{-1} (\varphi, (I - T_t) \varphi) 
= k(\varphi) 
< \infty
 \]
and $\varphi \one_\Omega \in D(k)$.
Then (\ref{eladty616;1}) follows  by taking the limit 
$t \downarrow 0$ in (\ref{eladty616;2}).

\ref{ladty616-2}$\Rightarrow$\ref{ladty616-1}.
Define the quadratic forms $k_\Omega$ on $L_2(\Omega)$ with form domain
\[
D(k_\Omega) = \{ \varphi \one_\Omega : \varphi \in D(k) \}
\]
and $k_\Omega(\varphi) = k(\varphi)$ for all $\varphi \in D(k_\Omega)$.
Define similarly the form $k_{\Omega^{\rm c}}$ on $L_2(\Omega^{\rm c})$.
Let $H_\Omega$ and $H_{\Omega^{\rm c}}$ be the associated self-adjoint operators.
Then it follows from (\ref{eladty616;1}) that 
$H = H_\Omega \oplus H_{\Omega^{\rm c}}$, where $H$ is the operator associated to $k$.
Then 
\[
T_t = e^{-t H_\Omega} \oplus e^{-t H_\Omega^{\rm c}}
\]
for all $t > 0$ and Statement~\ref{ladty616-1} follows.\hfill$\Box$

\ruimte

The method to prove separation is contained in the following lemma,
which  assumes the existence of suitable cut-off functions. 
We shall give several examples after the lemma.

\begin{lemma} \label{ladty608}
Let $\Omega$ be a measurable subset of $\Ri^d$.
Let $H_0$ be a viscosity operator with coefficients $c_{ij}$ on $L_2(\Ri^d)$.
Suppose there exist $\chi_1,\chi_2,\ldots \in L_{\infty;1}$ such that 
$0 \leq \chi_n \leq 1$ for all $n \in \Ni$, 
$\lim_{n \to \infty} \chi_n = \one_\Omega$ almost everywhere and 
\[
\lim_{n \to \infty} \int_W \sum_{i,j=1}^d c_{ij} \, (\partial_i \chi_n) \, (\partial_j \chi_n)
= 0
\]
for any compact subset $W$ of $\Ri^d$.
Then $S^{(0)}_t L_2(\Omega^{\rm c}) \subseteq L_2(\Omega^{\rm c})$ for all $t > 0$.
\end{lemma}
\proof\
It suffices to prove that
$S^{(0)}_t (L_2(\Omega^{\rm c}) \cap L_\infty(\Omega^{\rm c}))
          \subseteq L_2(\Omega^{\rm c})$ for all $t > 0$.
Let $\varphi \in L_2(\Omega^{\rm c}) \cap L_\infty(\Omega^{\rm c})$.
Define $X \colon \langle 0,\infty\rangle \to [0,\infty\rangle$ by
\[
X(t) = \int_\Omega |S^{(0)}_t \varphi|^2
= (S^{(0)}_t \varphi,\one_\Omega S^{(0)}_t \varphi)
\;\;\; .  \]
Then $\lim_{t \downarrow 0} X(t)=0$.
Moreover, $X$ is differentiable and
\[
X'(t) 
= - 2 (H_0 S^{(0)}_t \varphi,\one_\Omega S^{(0)}_t \varphi)
\]
for all $t > 0$.
Since $X \geq 0$ it suffices to show that $X'(t) \leq 0$ for all $t > 0$.

Fix $t > 0$ and $\tau \in C_c^\infty(\Ri^d)$ such that $\tau(x) = 1$ for all $x \in B_e(0\,;1)$ and 
$0 \leq \tau \leq 1$.
Define $\tau_R \in C_c^\infty(\Ri^d)$ for all $R > 0$   by 
$\tau_R(x) = \tau(R^{-1} x)$.
Set $\chi_{n,R} = \chi_n \, \tau_R$.
Then $\chi_{n,R} \in L_{2;1}(\Ri^d) \cap L_{\infty;1}(\Ri^d)$ for all $n \in \Ni$ and $R > 0$.

One has
\[
X'(t)= \lim_{R\to\infty} \lim_{n\to\infty} Y_{n,R}
\]
where
\[
Y_{n,R} = - 2  (H_0S^{(0)}_t\varphi , \chi_{n,R}^2 \, S^{(0)}_t\varphi)
\]
for all $n \in \Ni$ and $R \geq 1$.
Let $n \in \Ni$ and $R \geq 1$.
Since $H_\varepsilon$ converges in the strong resolvent sense to $H_0$
it follows that $\Psi(H_\varepsilon)$ converges strongly to $\Psi(H_0)$
for each bounded continuous function $\Psi$ on $\Ri$ by \cite{RS1}, Theorem~VIII.20(b).
Therefore
\[
Y_{n,R}
= \lim_{\varepsilon\to0} - 2 (H_\varepsilon S^{(\varepsilon)}_t\varphi , \chi_{n,R}^2 \,
                         S^{(\varepsilon)}_t\varphi)
\;\;\; .  \]
But integrating by parts  one finds
\begin{eqnarray*}
- 2 (H_\varepsilon S^{(\varepsilon)}_t\varphi , \chi_{n,R}^2 \,
                         S^{(\varepsilon)}_t\varphi)
& = & - 2 \sum_{i,j} ( (\partial_i S^{(\varepsilon)}_t \varphi),
   c^{(\varepsilon)}_{ij} \, \chi_{n,R}^2 \, (\partial_j S^{(\varepsilon)}_t\varphi))
\\*[5pt]
& & \hspace{10mm} {}
 - 4 \sum_{i,j} ((\partial_i S^{(\varepsilon)}_t\varphi),
              c^{(\varepsilon)}_{ij} \, \chi_{n,R} \, (\partial_j \chi_{n,R})
                                                       \,(S^{(\varepsilon)}_t\varphi))  \\[5pt]
&\leq& 2 \sum_{i,j} (S^{(\varepsilon)}_t\varphi,
           c^{(\varepsilon)}_{ij} \, (\partial_i \chi_{n,R}) \, (\partial_j \chi_{n,R})
                    \, S^{(\varepsilon)}_t\varphi)  
\;\;\; .  
\end{eqnarray*}
Since  $ \chi_{n,R} \in L_{\infty;1}$ one can again use strong resolvent convergence 
of the approximants to deduce that
\begin{eqnarray*}
Y_{n,R}
&\leq& 2 \sum_{i,j} (S^{(0)}_t\varphi, c_{ij} (\partial_i \chi_{n,R}) \, (\partial_j \chi_{n,R}) 
        \, S^{(0)}_t\varphi)  \\[5pt]
&\leq& 4 \sum_{i,j} (S^{(0)}_t\varphi, \tau_R^2 \, c_{ij} (\partial_i \chi_n) \, (\partial_j \chi_n) 
        \, S^{(0)}_t\varphi)   \\*[5pt]
& & \hspace{10mm} {}
   + 4 \sum_{i,j} (S^{(0)}_t\varphi, \chi_n^2 \, c_{ij} (\partial_i \tau_R) \, (\partial_j \tau_R) 
        \, S^{(0)}_t\varphi)  \\[5pt]
& \leq & 4 \|S^{(0)}_t\varphi\|_\infty^2 
         \int \sum_{i,j} \tau_R^2 \, c_{ij} (\partial_i \chi_n) \, (\partial_j \chi_n) 
   + a_1 \int_{B_e(0\,;R)^{\rm c}} |S^{(0)}_t\varphi|^2
\;\;\;,
\end{eqnarray*}
where $a_1 = 4 \|C\|_\infty \sum_{i=1}^d \|\partial_i \tau\|_\infty^2$
and we used  $R \geq 1$.
Since $\supp \tau_R$ is compact
one can use the assumption on the $\chi_n$ and 
take the limit $n \to \infty$.
Next take the limit $R \to \infty$.
Then  the last term also tends to zero.
Hence $X'(t) = \lim_{R\to\infty} \lim_{n\to\infty} Y_{n,R} \leq 0$ for all $t > 0$.\hfill$\Box$

\ruimte

The basic mechanism which leads to separation is one-dimensional.
Therefore we begin by  analyzing the one-dimensional situation in detail.
Subsequently we examine higher dimensions and describe different aspects that can occur.

In one-dimension there is one positive coefficient $c$. 
Suppose that $H_0$ is subelliptic.
Then $c>0$ almost everywhere,  $c^{-1/2}$ is integrable and 
$d_C$ is explicitly given by (\ref{eadt412}) by Theorem~\ref{tadty512}.
Now $d_C$ is finite-valued and 
the $L_2$ off-diagonal bounds of Proposition~\ref{padty311}
  are particularly simple. 
Then bounded open intervals are balls with a unique centre and radius 
and $\tilde d_C(I_1\,;I_2)=d_C(I_1\,;I_2)$ for each
pair of bounded open intervals $I_1$, $I_2$ because of the 
relation (\ref{eadt412}).
Hence Proposition~\ref{padty311} gives
\begin{equation}
|(\varphi_1,S^{(0)}_t\varphi_2)|\leq e^{-d_C(I_1;I_2)/(4t)}\|\varphi_1\|_2\,\|\varphi_2\|_2
\label{eadty5.100}
\end{equation}
for all $t>0$ and $\varphi_1\in L_2(I_1)$, $\varphi_2\in L_2(I_2)$ and all 
bounded open intervals $I_1$, $I_2$.
But a much stronger estimate is valid if the system separates into several disjoint subsystems.
The next proposition gives a criterion for the separation into two subsystems.
Recall that  
 an operator $S$ is positivity improving on $L_2$  if $S\varphi>0$ almost everywhere
whenever  $\varphi\geq0$ with $\varphi\neq0$.

\begin{prop}\label{padt402}
Let $x_0 \in \Ri$, $c\in L_\infty(\Ri)$ and $\alpha > 0$.
Assume  $c \geq 0$ almost everywhere, $c(x) > 0$ for all 
$x \in \langle x_0,x_0+\alpha]$, the function 
$c^{-1}$ is bounded on  $\langle x_0+\varepsilon,x_0+\alpha]$ for all 
$\varepsilon \in \langle0,\alpha\rangle$
and  that  $\int_{x_0}^{x_0+\alpha} c^{-1}=\infty$.
Then
\[
S^{(0)}_t L_2(-\infty,x_0) \subseteq L_2(-\infty,x_0)
\;\;\;\;\;{\rm and }
\;\;\;\;\;\;
S^{(0)}_t L_2(x_0, \infty) \subseteq L_2(x_0,\infty)
\]
for all $t>0$.
In particular the operator $S^{(0)}_t$ is not positivity improving for any $t>0$.
\end{prop}
\proof\
We may assume that $x_0 = 0$ and  $\alpha = 1$.
Secondly, for all $n\in\Ni$   define $\chi_n \colon \Ri \to [0,1]$ by
\[
\chi_n(x)
= \left\{ \begin{array}{ll}
   1 & \mbox{if } x\leq n^{-1} \;\;\; ,  \\[5pt]
   \eta_n^{-1} \, \eta(x) & \mbox{if } x\in\langle n^{-1},1\rangle \;\;\; , \\[5pt]
   0 & \mbox{if } x \geq 1 \;\;\; ,
          \end{array} \right.
\]
where
\[
\eta(x)=\int^1_x c^{-1}\;\;\;\;\;\;{\rm and}\;\;\;\;\;\eta_n=\eta(n^{-1})
\;\;\;.
\]
Note that $\chi_n$ is absolutely continuous and decreasing.
Moreover, $\lim_{n \to \infty} \chi_n = \one_{\langle-\infty,0]}$ pointwise
and $\chi_n' = - \eta_n^{-1} \, c^{-1} \, \one_{\langle n^{-1},1\rangle} \in L_\infty(\Ri)$ 
for all $n \in \Ni$.
Hence
\[
\|c\,(\chi_n')^2\|_1
= \int^1_{1/n} \eta_n^{-2} \, c^{-1} 
= \eta_n^{-1}
\]
for all $n \in \Ni$ and the proposition follows from Lemma~\ref{ladty608}.\hfill$\Box$

\ruimte

Of course the conclusion of the proposition is also valid if the conditions on $c$ 
on the right of $x_0$ are replaced by similar conditions on the left of $x_0$, i.e., 
if $c(x) > 0$ for all 
$x \in [x_0-\alpha,x_0\rangle$, the function 
$c^{-1}$ is bounded on the subsets $[x_0-\alpha,x_0-\varepsilon \rangle$ for all 
$\varepsilon \in \langle0,\alpha\rangle$
and  that  $\int_{x_0-\alpha}^{x_0} c^{-1}=\infty$.
Similar remarks are valid in other situations such as Propositions~\ref{padty611},
\ref{padty612} and \ref{padty5.5}.
We will not repeat this remark.

It is  clear from the Proposition~\ref{padt402}  and Corollary~\ref{cadty602}
that if the coefficient $c$  has several
zeros of the appropriate type then the system can split into several pieces.
For example one has the following.

\begin{prop} \label{padty611}
Let $x_1,x_2 \in \Ri$ with $x_1 < x_2$ and let $c \in L_\infty(\Ri)$ with 
$c \geq 0$ almost everywhere.
Suppose there exists an $\alpha > 0$ such that 
for each $k \in \{ 1,2 \} $ one has 
$c(x) > 0$ for all $x \in \langle x_k,x_k+\alpha]$
and the function 
$c^{-1}$ is bounded on $\langle x_k+\varepsilon,x_k+\alpha]$ for all 
$\varepsilon \in \langle0,\alpha\rangle$.
Further assume that  $\int_{x_k}^{x_k+\alpha} c^{-1}=\infty$.
Then $S^{(0)}$ leaves the subspaces $L_2(-\infty,x_1)$, $L_2(x_1,x_2)$ and $L_2(x_2,\infty)$ 
invariant and the semigroup is a direct sum of its restrictions to the subspaces.
Each such restriction is a positive contraction semigroup which extends to a 
contraction semigroup on each of the $L_p$-spaces and which is conservative
on the $L_\infty$-spaces.
\end{prop}

Since $L_2(-\infty,x_0)$ is the orthogonal component of $L_2(x_0,\infty)$ in $L_2(\Ri)$
one has a direct sum decomposition $L_2(\Ri)=L_2(-\infty,x_0)\oplus L_2(x_0,\infty)$.
Then Proposition~\ref{padt402} establishes that $S_t$ leaves the two subspaces invariant.
Hence one has a direct sum decomposition $S_t=S^{-(0)}_t\oplus S^{+(0)}_t$
where $S^{\pm(0)}$ denote the restrictions of $S^{(0)}$ to the appropriate subspaces. 
Then it follows straightforwardly that $H_0=H_{-0}\oplus H_{+0}$ where
 $H_{\pm0}$ denote the generators of  $S^{\pm(0)}$.
In particular it follows from these observations that the action of $S^{(0)}$ is not ergodic
and although the semigroup is positivity preserving it is not positivity improving.

The decomposition of $S^{(0)}$ implies that the $L_2$ off-diagonal bounds (\ref{eadty5.100}) 
can be strengthened since
\[
|(\varphi_1,S^{(0)}_t\varphi_2)|=0
\]
for all $t>0$ and $\varphi_1\in L_2(I_1)$, $\varphi_2\in L_2(I_2)$ with 
$I_1\subseteq \langle-\infty,x_0]$ and  $I_2\subseteq[x_0,\infty\rangle$.
This additional statement can be incorporated into (\ref{eadty5.100}) by replacing 
$d_C(I_1\,;I_2)$ by a set-theoretic `distance' by the method of Sturm 
\cite{Stu2}, page~237.
The definition of the  `distance' is superficially similar to the $d_C$-definition but
it is specifically adapted to the domain of the Dirichlet form $h_0$. 
In particular if separation takes place the `distance' between sets in different
components is infinity.

\smallskip

The separation phenomenon allows us to complete the discussion of Example~\ref{exadty4.1}
for $d=1$ and $\delta\in[1/2,1\rangle$.

\begin{exam}\label{exadty5.1}
Let $\delta \in \langle1/2,1\rangle$ and 
consider the one dimensional operator $H=-d\,c\,d$ with 
\[
c(x) = \Big( \frac{|x|^2}{1+|x|^2} \Big)^\delta
\;\;\; .  \]
We shall prove that $H_0$ is a subelliptic operator of order $1-\delta$.

Let $d$ denote the closed operator of differentiation on $L_2(\Ri)$ and
$d_\pm$  the corresponding   operators  on $L_2(\Ri_\pm)$ with 
domain $D(d_\pm) = \covsup{W}^{1,2}(\Ri_\pm)$, so with Dirichlet 
boundary conditions at the origin, where $\Ri_- = \langle-\infty,0\rangle$ and 
$\Ri_+ = \langle0,\infty\rangle$.
Then the adjoints ${d_\pm}^{\!\!*}$ are the operators of differentiation on $L_2(\Ri_\pm)$  with
domain $D({d_\pm}^{\!\!*}) = W^{1,2}(\Ri_\pm)$, so with no boundary condition.
Next define the form $h_{\pm,\varepsilon}$ on $L_2(\Ri_\pm)$ by
\[
h_{\pm,\varepsilon}(\varphi_\pm)=({d_\pm}^{\!\!*}\varphi_\pm,(c+\varepsilon){d_\pm}^{\!\!*}\varphi_\pm)
\]
and domain 
$D(h_{\pm,\varepsilon}) = D({d_\pm}^{\!\!*})$.
Moreover, set 
$h^{(N)}_\varepsilon$ on $L_2(\Ri)=L_2(\Ri_-)\oplus L_2(\Ri_+)$ by 
\[
h^{(N)}_\varepsilon(\varphi_-\oplus \varphi_+)
=h_{-,\varepsilon}(\varphi_-)+h_{+,\varepsilon}(\varphi_+)
\]
and domain $D(h^{(N)}_\varepsilon) = D(h_{-,\varepsilon}) \oplus D(h_{+,\varepsilon})$.
Now 
\[
h_\varepsilon(\varphi)=(d\varphi,(c+\varepsilon)d\varphi)
\]
for all $\varepsilon > 0$ with $D(h_\varepsilon)=D(d)$.
Then $h^{(N)}_\varepsilon\supseteq h_\varepsilon$ and 
$h_\varepsilon\geq h^{(N)}_\varepsilon$.
Thus $H_\varepsilon\geq  H_{-,\varepsilon}\oplus H_{+,\varepsilon}$ in the sense of quadratic forms
where $H_{\pm,\varepsilon}$ are the positive self-adjoint operators associated
with the forms $h_{\pm,\varepsilon}$.
Then by strong resolvent convergence one has 
$H_0\geq  H_{-,0}\oplus H_{+,0}$
where $H_{\pm,0}$ are the viscosity operators associated with the $H_{\pm,\varepsilon}$.
Therefore the problem of deriving subellipticity estimates on $H_0$ is reduced to 
deriving estimates on the operators $H_{\pm,0}$ on the subspaces $L_2(\Ri_\pm)$.

Let $l_\pm^D$ and $l_\pm^N$ be the forms on $L_2(\Ri_\pm)$ with domain
$D(l_\pm^D) = D(d_\pm)$ and $D(l_\pm^N) = D({d_\pm}^{\!\!*})$ given by
$l_\pm^D(\varphi) = \|d_\pm \varphi\|_2$ and $l_\pm^N(\varphi) = \|{d_\pm}^{\!\!*} \varphi\|_2$.
Then the self-adjoint operators $\Delta_{\pm D}$ and $\Delta_{\pm N}$ associated with 
$l_\pm ^D$ and $l_\pm ^N$ are called the Laplacians with Dirichlet and Neumann 
boundary conditions on $L_2(\Ri_\pm )$.
Then $\Delta_{\pm D}\geq(4x^2)^{-1}$ in the sense of quadratic forms.
This is the one-dimensional version of the estimate $\Delta\geq\sigma x^{-2}$
used in  higher dimensions in   Example~\ref{exadty4.1}.
It is a special case of the Hardy inequality (see, for example, \cite{Dav2}, Lemma~1.5.1).
Therefore $c\geq a_0\,(I + \Delta_{\pm D})^{-\delta}$, with $a_0>0$
and 
\[
h_{\pm ,\varepsilon}(\varphi)\geq 
a_0\,({d_\pm }^{\!\!*}\varphi,(I + \Delta_{\pm D})^{-\delta}{d_\pm }^{\!\!*}\varphi)
\]
for all $\varphi\in D({d_\pm }^{\!\!*})$.
But 
\[
{d_\pm }^{\!\!*}(I + \Delta_{\pm N})^{-\delta}\supseteq (I + \Delta_{\pm D})^{-\delta}{d_\pm }^{\!\!*}
\]
and so one has 
\[
h_{\pm ,\varepsilon}(\varphi)\geq  a_0\,({d_\pm }^{\!\!*}\varphi,{d_\pm }^{\!\!*}(I + \Delta_{\pm N})^{-\delta}\varphi)
=a_0\,(\Delta_{\pm N}^{1/2}\varphi,\Delta_{\pm N}^{1/2}(I + \Delta_{\pm N})^{-\delta}\varphi)
\]
for all $\varphi\in D({d_\pm }^{\!\!*})$.
Therefore
\[
H_{\pm ,\varepsilon}\geq a_0\,\Delta_{\pm N}(I + \Delta_{\pm N})^{-\delta}=
a_0\,\Delta_{\pm N}^{1-\delta}(\Delta_{\pm N} (I + \Delta_{\pm N})^{-1})^{\delta}
\]
in the sense of quadratic forms.
Then 
\[
H_{\pm ,0}\geq a_0\,\Delta_{\pm N}^{1-\delta}(\Delta_{\pm N} (I + \Delta_{\pm N})^{-1})^{\delta}
\]
by strong resolvent convergence.

Now let $\Delta_N = \Delta_{-N} \oplus \Delta_{+N}$.
We call $\Delta_N$ the Laplacian on $L_2(\Ri)$ with Neumann
boundary conditions at the origin.
Then 
\[
H_0\geq H_{-,0}\oplus H_{+,0}
\geq a_0\,\Delta_N^{1-\delta}(\Delta_N (I+\Delta_N)^{-1})^{\delta}
\]
and  by spectral theory there are $\mu,\nu>0$ such that 
\begin{equation}
H_0\geq \mu\,\Delta_N^{1-\delta}-\nu \,I
\;\;\; ,
\label{eadt415}
\end{equation}
i.e., the viscosity operator $H_0$ satisfies
the subellipticity condition relative to the Neumann Laplacian.
Up to now the estimates are in fact valid for all $\delta \in [0,1]$.

Set $\gamma = 1 - \delta$. 
Then $\gamma < 1/2$.
We strengthen the estimate (\ref{eadt415}) 
by use of standard properties of Sobolev spaces.
Since $\gamma < 1/2$ it follows by  straightforward argument involving $\sin$ and $\cos$
expansions 
that $D(\Delta_{\pm D}^{\gamma/2}) = D(\Delta_{\pm N}^{\gamma/2})$.
Hence there exists an $a_1 > 0$ such that 
$a_1 \, \|\Delta_{\pm D}^{\gamma/2} \varphi\|_2 \leq \|\Delta_{\pm N}^{\gamma/2} \varphi\|_2 + \|\varphi\|_2$
for all $\varphi \in D(\Delta_{\pm N}^{\gamma/2})$.
Then by dilation invariance one deduces that 
\[
a_1 \, \|\Delta_{\pm D}^{\gamma/2} \varphi\|_2 \leq \|\Delta_{\pm N}^{\gamma/2} \varphi\|_2
\]
for all $\varphi \in D(\Delta_{\pm N}^{\gamma/2})$.
Hence 
\begin{equation}
a_1^2 \, \Delta_{\pm D}^\gamma \leq \Delta_{\pm N}^\gamma
\label{eadtyx604;1}
\end{equation}
in the sense of quadratic forms.

Set $\Delta_D = \Delta_{-D} \oplus \Delta_{+D}$.
Then it follows from (\ref{eadtyx604;1}) that 
there is an $a_2 > 0$ such that 
$a_2 \, \Delta_D^\gamma \leq \Delta_N^\gamma$ in the sense of quadratic forms.
Recall that $\Delta$ denotes the Laplacian on $L_2(\Ri)$.
Then $\Delta \leq \Delta_D$ and therefore $\Delta^\gamma \leq \Delta_D^\gamma$.
Hence $a_2 \Delta^\gamma \leq \Delta_N^\gamma$.
Consequently  (\ref{eadt415}) gives 
\[
H_0\geq a_2 \, \mu \,\Delta^{1-\delta}-\nu I
\]
i.e., $H_0$ is  subelliptic  of order $1-\delta$.\hfill$\Box$
\end{exam}

Finally we consider the intermediate case $\delta = 1/2$.

\begin{exam}\label{xadty606}
Consider the one dimensional operator $H=-d\,c\,d$ with 
\[
c(x) = \Big( \frac{|x|^2}{1+|x|^2} \Big)^{1/2}
\;\;\; .  \]
We shall prove that $H_0$ is not a subelliptic operator of order $1/2$.

Since $c \in L_{\infty;1}(\Ri)$ the form $h$ is closable and $h_0 = \overline{h}$.
Let $\varphi \in C_c^\infty(\Ri)$ be such that $\varphi|_{[-1,1]} = 1$.
For all $n \in \Ni$ set $\varphi_n = \chi_n \varphi$, where $\chi_n$ is as
in the proof of Proposition~\ref{padt402}.
Then $\lim_{n \to \infty} \varphi_n = \varphi_-$ pointwise, where 
$\varphi_- = \varphi \one_{\langle-\infty,0]}$.
Then it follows from a modification of the calculations used to prove 
Proposition~\ref{padt402} that 
$\lim_{n \to \infty} h(\varphi_n) = \int_{-\infty}^0 c \, |\varphi'|^2$.

Next, $\varphi_-,\varphi_n \in L_1$ for all $n \in \Ni$ and 
$\lim_{n \to \infty} \varphi_n = \varphi_-$ in $L_1(\Ri)$.
Hence $\lim_{n \to \infty} \hat \varphi_n(p) = \hat \varphi_-(p)$ for all $p \in \Ri$.
Then
\[
\int_1^\infty dp \, |p| \, |\hat \varphi_-(p)|^2
\leq \liminf_{n \to \infty} \int_1^\infty dp \, |p| \, |\hat \varphi_n(p)|^2
\leq \liminf_{n \to \infty} \|\Delta^{1/4} \varphi_n\|_2^2
\;\;\; .  \]
Since $\varphi'$ vanishes in a neighbourhood of $0$ one calculates that for all $p \neq 0$
\[
\hat \varphi_-(p)
= i \, p^{-1} \, \varphi(0)
   - i \, p^{-1} ( \varphi'|_{\langle-\infty,0\rangle} ) \, \hat{\;}(p)
= i \, p^{-1} \, \varphi(0)
   - p^{-2} ( \varphi''|_{\langle-\infty,0\rangle} ) \, \hat{\;}(p)
\;\;\; .  \]
But $\varphi''|_{\langle-\infty,0\rangle} \in {\cal S}(\Ri)$ and 
therefore its Fourier transform
$( \varphi''|_{\langle-\infty,0\rangle} ) \, \hat{\;}$ is bounded.
Hence $\int_1^\infty dp \, |p| \, |\hat \varphi_-(p)|^2 = \infty$.
So $\liminf_{n \to \infty} \|\Delta^{1/4} \varphi_n\|_2^2 = \infty$.
Since $\sup_n h(\varphi_n) < \infty$ and $\|\varphi_n\|_2 \leq \|\varphi\|_2$ for 
all $n \in \Ni$ there are no $\mu,\nu > 0$ such that 
\[
h(\psi) \geq \mu \|\Delta^{1/4} \psi\|_2^2 - \nu \|\psi\|_2^2
\]
for all $\psi \in D(h) = L_{2;1} = D(\Delta^{1/2})$.
Thus $H_0$ is not  subelliptic  of order $1/2$.\hfill$\Box$
\end{exam}

One can now extend the conclusions of Example~\ref{exadty4.2} to $d=1$ and $\delta\in[1/2,1\rangle$.
For simplicity we only describe the case  of two zeros with the same order.

\begin{exam}\label{exadty5.2}
Let $\delta\in[1/2,1\rangle$ and $x_1<x_2$. 
Define $d \colon \Ri \to [0,\infty\rangle$ by
$d(x)=|x-x_1|\wedge |x-x_2|$.
Set 
\[
c(x)
= \Big( \frac{d(x)^2}{1+d(x)^2} \Big)^\delta
\;\;\; .  \]
Then separation takes place at both $x_1$ and $x_2$ by Proposition~\ref{padty611}.
Thus one has a decomposition of the semigroup  $S^{(0)}$ as a direct sum of its
components on the invariant subspaces $L_2(-\infty,x_1)$, 
$L_2(x_1,x_2)$ and $L_2(x_2,\infty)$.
Moreover, it follows by a slight variation of the arguments in Example~\ref{exadty5.1}
that one has subelliptic estimates
\[
H_0\geq\mu\,\Delta_N^{1-\delta}-\nu\,I
\]
where $\Delta_N$ is the direct sum of the Laplacians on the invariant subspaces with
Neumann boundary conditions at the endpoints $x_1$ and $x_2$.
Then, arguing as before, one has the full subelliptic estimate
\[
H_0\geq\mu'\,\Delta^{1-\delta}-\nu\,I
\;\;\;,
\]
with $\Delta$ the Laplacian on $L_2(\Ri)$, if $\delta\in\langle1/2,1\rangle$ but no such
estimate is valid if $\delta=1/2$.\hfill$\Box$
\end{exam}

The Examples~\ref{exadty4.1}, \ref{exadty5.1}, \ref{xadty606} and \ref{exadty5.2}
show that 
the properties of subelliptic operators are quite different
to those of strongly elliptic operators.
If formally $H = - d \, c \, d$ and $H_0$ is the viscosity operator
then $H$ is strongly elliptic if $c \in C^2_b(\Ri)$ by Proposition~\ref{padty513}.
If $\delta \in [1/2,1\rangle$ and 
\[
c_\delta(x) = \Big( \frac{|x|^2}{1+|x|^2} \Big)^\delta
\]
then $c_\delta$ is differentiable and $c_\delta' \in C^{2\delta-1}(\Ri)$,
with obvious modifications if $\delta = 1/2$.
If $H_0$ is the viscosity operator with coefficient $c_\delta$ then 
$H_0$ is subelliptic by Example~\ref{exadty5.1} if $\delta > 1/2$.
Hence the viscosity semigroup $S^{(0)}$ has a bounded kernel $K^{(0)}$ by
Proposition~\ref{padty4.1}.
But separation takes place by Proposition~\ref{padt402}.
Therefore
\begin{itemize}
\item
for all $t > 0$ the operator $S^{(0)}_t$ is not positivity improving
(Proposition~\ref{padt402})
\end{itemize}
and in particular
\begin{itemize}
\item
for all $t > 0$ the  kernel $K^{(0)}_t$ is not strictly positive.
\end{itemize}
On the other hand, for strongly elliptic operators with real measurable bounded 
coefficients the kernel is always strictly positive and H\"older continuous 
by Nash and De Giorgi.
Now, however, 
\begin{itemize}
\item
for all $t > 0$ the  kernel $K^{(0)}_t$ is not continuous.
\end{itemize}

Separation implies that 
$K^{(0)}_t(-x\,;y) = 0$ whenever $x,y > 0$.
Hence if $K^{(0)}_t$ is continuous then  $K^{(0)}_t(0\,;0) = \lim_{x \downarrow 0} K^{(0)}_t(-x\,;x) = 0$, in 
contradiction to Corollary~\ref{ctdy3.3}.

The example also shows that Fefferman and Phong's geometric criterion for subellipticity,
\cite{FP} Theorem~1, is sensitive to smoothness of the coefficients.
Theorem~1 in \cite{FP} states that if $c_{ij} \in C^\infty(\Ri^d)$ with $C \geq 0$ everywhere and if 
$\delta \in \langle0,1]$ then $H$ is subelliptic of order $\delta$ if and only if
\begin{equation}
\mbox{for all } x \in \Ri^d \mbox{ there is an } a > 0 \mbox{ such that }
d_C(x\,;y) \leq a \, |x-y|^\delta \mbox{ for all } y \in B_e(x\,;1)
\label{eSadty6;1}
\end{equation}
In Theorem~\ref{tadty512} we showed in one-dimension, for operators with measurable 
coefficients, that 
subellipticity of order~$\delta$ implies (\ref{eSadty6;1}) and, in addition, that one can choose
$a$ independent of~$x$.
The converse is not valid in general since
\begin{itemize}
\item
if $\delta = 1/2$ then (\ref{eSadty6;1}) is valid with $\delta = 1/2$, but $H_0$ is 
not subelliptic of order $1/2$
\end{itemize}
by Example~\ref{xadty606}.

This shows that the claim in \cite{Stu2}, Theorem~4.3.(ii), requires additional assumptions.

Kusuoka and Stroock \cite{KuS} proved under suitable conditions that 
the kernel satisfies large time Euclidean Gaussian bounds involving the dimension $d$
of the underlying space $\Ri^d$, 
But if $\delta \in \langle1/2,1\rangle$ and $H_0$ is the operator in 
Example~\ref{exadty5.2} then $H_0$ is subelliptic but
\begin{itemize}
\item
there are no $a,b > 0$ such that 
$K^{(0)}_t(x\,;y) \leq a \, t^{-d/2} e^{-b |x-y|^2 t^{-1}}$
for all $x,y \in \Ri^d$ and $t \geq 1$.
\end{itemize}
In fact
\begin{itemize}
\item
there is no $a > 0$ such that 
$\|K^{(0)}_t\|_\infty \leq a \, t^{-d/2}$
for all $t \geq 1$.
\end{itemize}
In the notation of Example~\ref{exadty5.2} one has $S^{(0)}_t \one_{[x_1,x_2]} = \one_{[x_1,x_2]}$.
Hence 
\[
1 = \|S^{(0)}_t \one_{[x_1,x_2]}\|_\infty
\leq \|S^{(0)}_t\|_{1 \to \infty} \|\one_{[x_1,x_2]}\|_1
= \|K^{(0)}_t\|_\infty  \|\one_{[x_1,x_2]}\|_1
\]
for all $t > 0$ and $\|K^{(0)}_t\|_\infty \geq (x_2 - x_1)^{-1}$.
This shows that the subellipticity condition in Theorem~\ref{tadty411} is not sufficient
to prove (\ref{etdy4.111}) and something more is necessary.

It is also evident that the kernel cannot satisfy Gaussian lower bounds because 
the kernel is not strictly positive.

\smallskip

Although  the foregoing description of separation is restricted to one-dimension 
one can use the  mechanism of Proposition~\ref{padt402} to  construct examples 
in higher dimensions for which one obtains separation.
As a first illustration we consider the separation of a compact subset which we 
take to be the Euclidean ball $B_e(0\,;1)=\{x\in\Ri^d: |x|<1\}$.

\begin{prop} \label{padty612}
Let $c \in L_\infty([0,\infty\rangle)$ with $c \geq 0$ almost everywhere.
Assume $c(x) > 0$ for all 
$x \in \langle 1,1+\alpha]$ and the function 
$c^{-1}$ is bounded on  $\langle 1+\varepsilon,1+\alpha]$ for all 
$\varepsilon \in \langle0,\alpha\rangle$
and  that  $\int_1^{1+\alpha} c^{-1}=\infty$.
Define $c_{ij} \colon \Ri^d \to \Ri$ by
$c_{ij}(x) = c(|x|) \, \delta_{ij}$.
Let $H_0$ be the viscosity operator with coefficients $c_{ij}$.
Then $S^{(0)}_t L_2(B_e(0\,;1)) \subseteq L_2(B_e(0\,;1))$ for all $t > 0$.
\end{prop}
\proof\ 
We may assume that $\int_1^2 c^{-1} = \infty$.
For all $n \in \Ni$ define $\chi_n \colon \Ri^d \to [0,1]$ by
\[
\chi_n(x)
= \left\{ \begin{array}{ll}
   1 & \mbox{if } |x| \leq 1 + n^{-1} \;\;\; ,  \\[5pt]
   \eta_n^{-1} \, \eta(|x|) & \mbox{if } |x| \in \langle 1+ n^{-1},2\rangle \;\;\; , \\[5pt]
   0 & \mbox{if } |x| \geq 2 \;\;\; ,
          \end{array} \right.
\]
where
\[
\eta(r)=\int^2_r c^{-1}\;\;\;\;\;\;{\rm and}\;\;\;\;\;\eta_n=\eta(1+n^{-1})
\;\;\;.
\]
Then $\lim_{n \to \infty} \chi_n = \one_{\overline{B_e(0\,;1)}}$ pointwise
and $\chi_n \in L_{\infty;1}$ for all $n$.
Next, 
\begin{eqnarray*}
\sum^d_{i,j=1} c_{ij}(x) \, (\partial_i\chi_n)(x) \, (\partial_j\chi_n)(x)
& = & \eta_n^{-2} \, c(|x|)^{-1}  \one_{\langle1+n^{-1},2\rangle}(|x|)
\end{eqnarray*}
for almost every $x \in \Ri^d$.
Therefore 
\[
\int_{\Ri^d} \sum_{i,j=1}^d c_{ij} \, (\partial_i \chi_n) \, (\partial_j \chi_n)
= |B_e(0\,;1)| \, \eta_n^{-2} \, \int_{1+n^{-1}}^2 dr \, r^{d-1} \, c(r)^{-1}
\leq |B_e(0\,;1)| \, \eta_n^{-1}
\]
for all $n \in \Ni$ and the proposition follows from Lemma~\ref{ladty608}.\hfill$\Box$

\ruimte

One can use the reasoning in Examples~\ref{exadty4.1} and \ref{exadty4.2} to obtain
explicit examples of this form of separation.

\begin{exam}\label{exadty5.3}
Define $d \colon \Ri^d \to \Ri$ by  $d(x)= \min \{ |x-y| : |y| =  1 \}$.
Let $\delta\in[1/2,1\rangle$ .
Assume $(c_{ij})\geq ( c\,\delta_{ij})$ 
where 
\[
c(x) = \Big( \frac{d(x)^2}{1+d(x)^2} \Big)^\delta
\;\;\; .  \]
Then separation takes place and the system
factors into a component in the unit ball $B_e(0\,;1)$ and a component in its
complement.
The subellipticity properties of $H_0$ are similar to the one-dimensional examples.
The operator is subelliptic with respect to the operator formed as the direct sum
of the Laplacian in the interior of the ball, with Neumann boundary conditions 
on the boundary, and the Laplacian in the exterior, with the same boundary conditions,
by \cite{Dav2} Theorem~1.5.4.
Moreover, the operator is subelliptic with respect to the full Laplacian on $\Ri^d$ 
if $\delta\in\langle1/2,1\rangle$ but not if $\delta=1/2$.

Since $\one_{\overline{B_e(0\,;1)}} \in L_1(\Ri^d)$ one deduces again that 
$\|K^{(0)}_t\|_\infty \geq |B_e(0\,;1)|^{-1}$ for all $t > 0$.
Therefore it is   not possible to have  bounds $\|K^{(0)}_t\|_\infty \leq a \, t^{-d/2}$ 
uniformly for all $t \geq 1$.\hfill$\Box$
\end{exam}

One can also construct examples in higher dimensions where the system separates across an
infinitely extended surface.
Note that we assume no bounds on the derivatives of the function $\Phi$ in the next proposition.

\begin{prop}\label{padty5.5}
Let $\Phi\in C^1(\Ri^{d-1})$.
Define
\[
\Omega=\{ (y,z) \in \Ri^{d-1} \times \Ri : z \leq \Phi(y) \}
\]
Let  $c\in L_\infty(\Ri)$ with $c \geq 0$ almost everywhere and $\alpha > 0$.
Assume  $c(x) > 0$ for all 
$x \in \langle 0,\alpha]$ and the function 
$c^{-1}$ is bounded on the subsets $\langle\varepsilon,\alpha]$ for all 
$\varepsilon \in \langle0,\alpha\rangle$
and  that  $\int_0^{\alpha} c^{-1}=\infty$.
Define $c_{ij} \colon \Ri^{d-1} \times \Ri \to \Ri$ by
$c_{ij}(y,z) = c(z - \Phi(y)) \, \delta_{ij}$.
Then 
\[
S^{(0)}_t L_2(\Omega) \subseteq L_2(\Omega)
\]
for all $t > 0$.
\end{prop}
\proof\
We may assume that $\alpha=1$.
Let  $\chi_n$ denote the functions introduced in the proof of Proposition~\ref{padt402}
and define $\tilde\chi_n\in L_{\infty;1}(\Ri^{d-1} \times \Ri)$ by 
$\tilde\chi_n(y,z) = \chi_n(z - \Phi(y))$.
Then $\lim_{n \to \infty} \|c \, (\chi_n')^2\|_1 = 0$.
Moreover, $\lim_{n \to \infty} \tilde\chi_n = \one_\Omega$ pointwise.
For almost every $(y,z) \in \Ri^{d-1} \times \Ri$ one has
\[
\sum_{i,j=1}^d c_{ij}(y,z) \, (\partial_i \chi_n)(y,z) \, (\partial_j \chi_n)(y,z) 
= c(z - \Phi(y)) \, |\chi_n'(z - \Phi(y))|^2 \, (1 + |(\nabla \Phi)(y)|^2)
\;\;\; .  \]
Since $\lim_{n \to \infty} \|c \, (\chi_n')^2\|_1 = 0$
it follows that
\[
\lim_{n \to \infty} \int_W \sum_{i,j=1}^d c_{ij} \, (\partial_i \chi_n) \, (\partial_j \chi_n)
= 0
\]
for any compact subset $W$ of $\Ri^d$.
Then the proposition is a consequence of  Lemmas~\ref{ladty608} and \ref{ladty601}.\hfill$\Box$

\ruimte

Although Proposition~\ref{padty5.5} only deals with separation by one surface
one may easily extend the reasoning  to separation across several surfaces.
An interesting situation occurs if one has two disjoint surfaces such as 
$z = \pm  \Phi(y)$
where $\Phi > 0$ and $\lim_{|y| \to \infty} \Phi(y) = 0$.
Then the system splits into three components 
$\Omega_\pm= \{ (y,z) \in \Ri^{d-1} \times \Ri : \pm z > \Phi(y) \} $
and $\Omega_0 = \{ (y,z) \in \Ri^{d-1} \times \Ri : |z| \leq \Phi(y) \} $.
All three components  are unbounded but the  $\Omega_0$ component
can have finite volume and the corresponding operator $H_0|_{\Omega_0}$
can have compact resolvent. 
Then again $\|K^{(0)}_t\|_\infty\geq |\Omega_0|^{-1}$ and one cannot have 
Euclidean Gaussian bounds.
In this situation it is also not possible to have a small 
$t$ power behaviour of the semigroup.

\begin{exam}\label{exadty5.4}
Let $\Phi \in C^1(\Ri^{d-1})$ with $\Phi > 0$ and $\lim_{|y| \to \infty} \Phi(y) = 0$.
Set $\Omega_\pm= \{ (y,z) \in \Ri^{d-1} \times \Ri : \pm z > \Phi(y) \} $
and $\Omega_0 = \{ (y,z) \in \Ri^{d-1} \times \Ri : |z| \leq \Phi(y) \} $.
Let $\delta \in [1/2,1\rangle$.
Define $c_{ij} \colon \Ri^{d-1} \times \Ri \to \Ri$ by
\[
c_{ij}(y,z) = 
\Big( \frac{ (|z| - \Phi(y))^2 }{1 + (|z| - \Phi(y))^2 } \Big)^\delta \, \delta_{ij}
\;\;\; .  \]
Then it follows as in the proof of Proposition~\ref{padty5.5}
and by Corollary~\ref{cadty602} that
one has separation into the three components $\Omega_\pm$ and 
$\Omega_0$.
Now $S^{(0)}$ leaves $L_2(\Omega_0)$ invariant and $S^{(0)} \one_{\Omega_0} = \one_{\Omega_0}$
for all $t > 0$.
Let $a, \alpha>0$ and $\omega\geq0$.
Suppose 
\begin{equation}
\|S^{(0)}_t\|_{1\to2}\leq a\,t^{-\alpha}e^{\omega t}
\label{eadty6.137}
\end{equation}
for all $t>0$.
Fix $m \in \Ni$ with $m>\alpha$.
Then there are $a_1,R>0$ such that 
\[
\|(I+r^2H_0)^{-m}\|_{1\to2}\leq a_1\,r^{-2\alpha}
\]
for all $r\in\langle0,R]$.
(See the proof of Theorem~\ref{tadt302}.)
Therefore the estimate (\ref{etadt302;1}) gives 
\[
|B_C((y,z)\,;r)\cap\Omega_0|
\geq a_1^{-1}\,r^{2\alpha}
\]
for $(y,z) \in \Omega_0$ and all $r\in\langle0,R]$.
But  with $r$ fixed it follows that  
\[
\lim_{|y| \to \infty} |B_C((y,0)\,;r)\cap\Omega_0| = 0
\;\;\; .  \]
This is a contradiction.
Therefore the bounds (\ref{eadty6.137}) are not possible. 
In particular $H_0$ cannot be subelliptic.\hfill$\Box$
\end{exam}

Note that the last example differs in character from the earlier ones
insofar one has  $c_{ij}(y,\Phi(y)) = 0$ for all $y \in \Ri^{d-1}$,
i.e., the coefficients are zero on an unbounded hypersurface.

\section*{Acknowledgements}

This work was supported by the Australian Research Council's Discovery
Grant Program.
Parts of the work were carried out whilst the first author 
was visiting the Centre of Mathematics and its Applications at the ANU, with ARC support, 
and  the second author was visiting the Eindhoven  University of Technology.
The third and fourth authors participated in the collaboration with the support of the ARC
whilst on leave at the ANU .

All four  authors thank Barry Simon and Alan McIntosh for helpful information regarding
the  low-dimensional examples in Sections~\ref{Sadty4} and \ref{Sadty5}.


\begin{thebibliography}{FLW95} \addcontentsline{toc}{chapter}{Bibliography}
  \markboth{Bibliography}{Bibliography}

\bibitem[Aro]{Aro}
{\sc Aronson, D.G.}, Bounds for the fundamental solution of a parabolic
  equation.
\newblock {\em Bull.\ Amer. Math.\ Soc.} {\bf 73} (1967),  890--896.

\bibitem[Aus1]{Aus1}
{\sc Auscher, P.}, Regularity theorems and heat kernels for elliptic operators.
\newblock {\em J. London Math.\ Soc.} {\bf 54} (1996),  284--296.

\bibitem[Aus2]{Aus2}
\leavevmode\vrule height 2pt depth -1.6pt width 23pt, On necessary and
  sufficient conditions for $L^p$-estimates of Riesz transforms associated to
  elliptic operators on $\Ri^n$ and related estimates.
\newblock Research report, Preprint. Univ.\ de Paris-Sud, 2004.

\bibitem[BeD]{BeD}
{\sc Beurling, A. {\rm and} Deny, J.}, Espaces de Dirichlet. I. Le cas
  \'el\'ementaire.
\newblock {\em Acta Math.} {\bf 99} (1958),  203--224.

\bibitem[BiM]{BM}
{\sc Biroli, M. {\rm and} Mosco, U.}, A Saint-Venant type principle for
  Dirichlet forms on discontinuous media.
\newblock {\em Ann.\ Mat.\ Pura Appl.} {\bf 169} (1995).

\bibitem[BoH]{BH}
{\sc Bouleau, N. {\rm and} Hirsch, F.}, {\em Dirichlet forms and analysis on
  Wiener space}, vol.\ 14 of de Gruyter Studies in Mathematics.
\newblock Walter de Gruyter \& Co., Berlin, 1991.

\bibitem[Bra]{Bra}
{\sc Braides, A.}, {\em $\Gamma$-convergence for beginners}, vol.\ 22 of Oxford
  Lecture Series in Mathematics and its Applications.
\newblock Oxford University Press, Oxford, 2002.

\bibitem[CKS]{CKS}
{\sc Carlen, E.A., Kusuoka, S. {\rm and} Stroock, D.W.}, Upper bounds for
  symmetric Markov transition functions.
\newblock {\em Ann.\ Inst.\ Henri Poincar\'e} {\bf 23} (1987),  245--287.

\bibitem[CGT]{CGT}
{\sc Cheeger, J., Gromov, M. {\rm and} Taylor, M.}, Finite propagation speed,
  kernel estimates for functions of the Laplace operator, and the geometry of
  complete Riemannian manifolds.
\newblock {\em J. Differential Geom.} {\bf 17} (1982),  15--53.

\bibitem[Cou]{Cou4}
{\sc Coulhon, T.}, Off-diagonal heat kernel lower bounds without Poincar\'e.
\newblock {\em J. London Math.\ Soc.} {\bf 68} (2003),  795--816.

\bibitem[Dav1]{Dav7}
{\sc Davies, E.B.}, Explicit constants for Gaussian upper bounds on heat
  kernels.
\newblock {\em Amer.\ J. Math.} {\bf 109} (1987),  319--333.

\bibitem[Dav2]{Dav2}
\leavevmode\vrule height 2pt depth -1.6pt width 23pt, {\em Heat kernels and
  spectral theory}.
\newblock Cambridge Tracts in Mathematics 92. Cambridge University Press,
  Cambridge etc., 1989.

\bibitem[Dav3]{Dav12}
\leavevmode\vrule height 2pt depth -1.6pt width 23pt, Heat kernel bounds,
  conservation of probability and the Feller property.
\newblock {\em J. Anal.\ Math.} {\bf 58} (1992),  99--119.
\newblock Festschrift on the occasion of the 70th birthday of Shmuel Agmon.

\bibitem[DeG]{DG}
{\sc De Giorgi, E.}, Sulla differenziabilit\`a e l'analiticit\`a delle
  estremali degli integrali multipli regolari.
\newblock {\em Mem.\ Accad.\ Sci.\ Torino cl.\ Sci.\ Fis.\ Mat.\ Nat.} {\bf 3}
  (1957),  25--43.

\bibitem[EkT]{ET}
{\sc Ekeland, I. {\rm and} Temam, R.}, {\em Convex analysis and variational
  problems}.
\newblock North-Holland Publishing Co., Amsterdam, 1976.

\bibitem[ElR1]{ER22}
{\sc Elst, A.F.M. ter {\rm and} Robinson, D.W.}, Local lower bounds on heat
  kernels.
\newblock {\em Positivity} {\bf 2} (1998),  123--151.

\bibitem[ElR2]{ER15}
\leavevmode\vrule height 2pt depth -1.6pt width 23pt, Second-order subelliptic
  operators on Lie groups I: complex uniformly continuous principal
  coefficients.
\newblock {\em Acta Appl.\ Math.} {\bf 59} (1999),  299--331.

\bibitem[ElR3]{ER18}
\leavevmode\vrule height 2pt depth -1.6pt width 23pt, Second-order subelliptic
  operators on Lie groups II: real measurable principal coefficients.
\newblock In {\sc Balakrishnan, A.V.}, ed., {\em Proceedings for the First
  International Conference of Semigroups of Operators: Theory and Applications,
  Newport Beach, California}, vol.\ 42 of Progress in nonlinear differential
  equations and their applications. Birkh{\"a}user Verlag, Basel, 2000,
  103--124.

\bibitem[ElR4]{ER27}
\leavevmode\vrule height 2pt depth -1.6pt width 23pt, Global estimates for
  H{\"o}rmander sublaplacians, 2005.
\newblock In preparation.

\bibitem[ERZ]{ERZ1}
{\sc Elst, A.F.M. ter, Robinson, D.W., {\rm and} Zhu, Y.}, Positivity and
  ellipticity.
\newblock {\em Proc.\ Amer.\ Math.\ Soc.} (2005).
\newblock To appear.

\bibitem[FKS]{FKS}
{\sc Fabes, E.B., Kenig, C.E. {\rm and} Serapioni, R.P.}, The local
  regularity of solutions of degenerate elliptic equations.
\newblock {\em Comm.\ Part.\ Diff.\ Eq.} {\bf 7} (1982),  77--116.

\bibitem[FaS]{FaS}
{\sc Fabes, E.B. {\rm and} Stroock, D.W.}, A new proof of Moser's parabolic
  Harnack inequality using the old ideas of Nash.
\newblock {\em Arch.\ Rat.\ Mech.\ and Anal.} {\bf 96} (1986),  327--338.

\bibitem[FeP]{FP}
{\sc Fefferman, C. {\rm and} Phong, D.H.}, Subelliptic eigenvalue problems.
\newblock In {\em Conference on harmonic analysis in honor of Antoni Zygmund},
  Wadsworth Math.\ Ser.,  590--606. Wadsworth, Belmont, CA, 1983.

\bibitem[FeS]{FSC}
{\sc Fefferman, C.L. {\rm and} S{\'a}nchez-Calle, A.}, Fundamental solutions
  for second order subelliptic operators.
\newblock {\em Ann.\ Math.} {\bf 124} (1986),  247--272.

\bibitem[Fra]{Fra}
{\sc Franchi, B.}, Weighted Sobolev-Poincar\'e inequalities and pointwise
  estimates for a class of degenerate elliptic equations.
\newblock {\em Trans.\ Amer.\ Math.\ Soc.} {\bf 327} (1991),  125--158.

\bibitem[FrL]{FL}
{\sc Franchi, B. {\rm and} Lanconelli, E.}, H{\"o}lder regularity theorem for a
  class of linear nonuniformly elliptic operators with measurable coefficients.
\newblock {\em Ann.\ Scuola Norm.\ Sup.\ Pisa Cl.\ Sci.} {\bf 10} (1983),
  523--541.

\bibitem[FLW]{FLW}
{\sc Franchi, B., Lu, G. {\rm and} Wheeden, R.L.}, Weighted Poincar\'e
  inequalities for H{\"o}rmander vector fields and local regularity for a class
  of degenerate elliptic equations.
\newblock {\em Potential Anal.} {\bf 4} (1995),  361--375.

\bibitem[Fri]{Fri}
{\sc Friedman, A.}, {\em Partial differential equations of parabolic type}.
\newblock Prentice-Hall, Inc., Englewood Cliffs, N.J., 1964.

\bibitem[FOT]{FOT}
{\sc Fukushima, M., Oshima, Y. {\rm and} Takeda, M.}, {\em Dirichlet forms and
  symmetric Markov processes}, vol.\ 19 of de Gruyter Studies in Mathematics.
\newblock Walter de Gruyter \& Co., Berlin, 1994.

\bibitem[Gaf]{Gaf}
{\sc Gaffney, M.P.}, The conservation property of the heat equation on
  Riemannian manifolds.
\newblock {\em Comm.\ Pure Appl.\ Math.} {\bf 12} (1959),  1--11.

\bibitem[GeS]{GS}
{\sc Gel'fand, I.M. {\rm and} Shilov, G.E.}, {\em Generalized functions},
  vol.\ 2.
\newblock Academic Press, New York, 1968.

\bibitem[Gia1]{Gia1}
{\sc Giaquinta, M.}, {\em Multiple integrals in the calculus of variations and
  nonlinear elliptic systems}.
\newblock Annals of Mathematics Studies 105. Princeton University Press,
  Princeton, 1983.

\bibitem[Gia2]{Gia}
\leavevmode\vrule height 2pt depth -1.6pt width 23pt, {\em Introduction to
  regularity theory for nonlinear elliptic systems}.
\newblock Lectures in Mathematics ETH Z{\"u}rich. Birkh{\"a}user Verlag, Basel
  etc., 1993.

\bibitem[GiT]{GT}
{\sc Gilbarg, D. {\rm and} Trudinger, N.S.}, {\em Elliptic partial
  differential equations of second order}.
\newblock Second edition, Grundlehren der mathematischen Wissenschaften 224.
  Springer-Verlag, Berlin etc., 1983.

\bibitem[Gri]{Gri3}
{\sc Grigor'yan, A.}, Estimates of heat kernels on Riemannian manifolds.
\newblock In {\em Spectral theory and geometry $($Edinburgh, 1998$)$}, vol.\ 273 of
  London Math.\ Soc.\ Lecture Note Ser.,  140--225. Cambridge Univ.\ Press,
  Cambridge, 1999.

\bibitem[H{\"o}r1]{Hor1}
{\sc H{\"o}rmander, L.}, Hypoelliptic second order differential equations.
\newblock {\em Acta Math.} {\bf 119} (1967),  147--171.

\bibitem[H{\"o}r2]{Hor4}
\leavevmode\vrule height 2pt depth -1.6pt width 23pt, {\em The analysis of
  linear partial differential operators. I{}. Distribution theory and Fourier
  analysis}, vol.\ 256 of Grundlehren der Mathematischen Wissenschaften.
\newblock Springer-Verlag, Berlin, 1983.

\bibitem[H{\"o}r3]{Hor5}
\leavevmode\vrule height 2pt depth -1.6pt width 23pt, {\em The analysis of
  linear partial differential operators. II{}. Differential operators with
  constant coefficients}, vol.\ 257 of Grundlehren der Mathematischen
  Wissenschaften.
\newblock Springer-Verlag, Berlin, 1983.

\bibitem[H{\"o}r4]{Hor6}
\leavevmode\vrule height 2pt depth -1.6pt width 23pt, {\em The analysis of
  linear partial differential operators. III{}. Pseudodifferential operators},
  vol.\ 274 of Grundlehren der Mathematischen Wissenschaften.
\newblock Springer-Verlag, Berlin, 1985.

\bibitem[H{\"o}r5]{Hor7}
\leavevmode\vrule height 2pt depth -1.6pt width 23pt, {\em The analysis of
  linear partial differential operators. IV{}. Fourier integral operators},
  vol.\ 275 of Grundlehren der Mathematischen Wissenschaften.
\newblock Springer-Verlag, Berlin, 1985.

\bibitem[JeS]{JSC1}
{\sc Jerison, D. {\rm and} S{\'a}nchez-Calle, A.}, Subelliptic, second order
  differential operators.
\newblock In {\sc Berenstein, C.A.}, ed., {\em Complex analysis III}, Lecture
  Notes in Mathematics 1277. Springer-Verlag, Berlin etc., 1987,  46--77.

\bibitem[Jos]{Jos}
{\sc Jost, J.}, Nonlinear Dirichlet forms.
\newblock In {\em New directions in Dirichlet forms}, vol.\ 8 of AMS/IP Stud.\
  Adv.\ Math.,  1--47. Amer.\ Math.\ Soc., Providence, RI, 1998.

\bibitem[Kat]{Kat1}
{\sc Kato, T.}, {\em Perturbation theory for linear operators}.
\newblock Second edition, Grundlehren der mathematischen Wissenschaften 132.
  Springer-Verlag, Berlin etc., 1984.

\bibitem[KuS]{KuS}
{\sc Kusuoka, S. {\rm and} Stroock, D.}, Long time estimates for the heat
  kernel associated with a uniformly subelliptic symmetric second order
  operator.
\newblock {\em Ann.\ Math.} {\bf 127} (1988),  165--189.

\bibitem[LaM]{LaM}
{\sc Lanconelli, E. {\rm and} Morbidelli, D.}, On the Poincar\'e inequality for
  vector fields.
\newblock {\em Ark.\ Mat.} {\bf 38} (2000),  327--342.

\bibitem[MaR]{MR}
{\sc Ma, Z.M. {\rm and} R{\"o}ckner, M.}, {\em Introduction to the theory of
  (non symmetric) Dirichlet Forms}.
\newblock Universitext. Springer-Verlag, Berlin etc., 1992.

\bibitem[Mar]{Mark}
{\sc Markushevich, A.I.}, {\em Theory of functions of a complex variable.
  Vols.\ I, II, III}.
\newblock Chelsea Publishing Co., New York, 1977.

\bibitem[Mas]{DalM}
{\sc Maso, G.D.}, {\em An introduction to {$\Gamma$}-convergence}, vol.\ 8 of
  Progress in Nonlinear Differential Equations and their Applications.
\newblock Birkh{\"a}user Boston Inc., Boston, MA, 1993.

\bibitem[Mel]{Melr}
{\sc Melrose, R.}, Propagation for the wave group of a positive subelliptic
  second-order differential operator.
\newblock In {\em Hyperbolic equations and related topics (Katata/Kyoto,
  1984)},  181--192. Academic Press, Boston, MA, 1986.

\bibitem[Mos]{Mosco}
{\sc Mosco, U.}, Composite media and asymptotic Dirichlet forms.
\newblock {\em J. Funct.\ Anal.} {\bf 123} (1994),  368--421.

\bibitem[Mus]{Mus}
{\sc Mustapha, S.}, Sous ellipticit\'e dans le cadre du calcul $S(m,g)$.
\newblock {\em Comm.\ Part.\ Diff.\ Eq.} {\bf 19} (1994),  245--275.

\bibitem[MuV]{MV1}
{\sc Mustapha, S. {\rm and} Varopoulos, N.T.}, Comparaison H{\"o}lderienne des
  distances sous-elliptiques et calcul $S(m,g)$.
\newblock {\em Potential Anal.} {\bf 4} (1995),  415--428.

\bibitem[Nas]{Nash}
{\sc Nash, J.}, Continuity of solutions of parabolic and elliptic equations.
\newblock {\em Amer.\ J. Math.} {\bf 80} (1958),  931--954.

\bibitem[OlR]{OR}
{\sc Ole{\u{\i}}nik, O.A. {\rm and} Radkevi\v{c}, E.V.}, {\em Second order
  equations with nonnegative characteristic form}.
\newblock Amer.\ Math.\ Soc., Providence, 1973.

\bibitem[Rob]{Robm}
{\sc Robinson, D.W.}, {\em Elliptic operators and Lie groups}.
\newblock Oxford Mathematical Monographs. Oxford University Press, Oxford etc.,
  1991.

\bibitem[ReS1]{RS1}
{\sc Reed, M. {\rm and} Simon, B.}, {\em Methods of modern mathematical physics
  I. Functional analysis}.
\newblock Academic Press, New York etc., 1972.

\bibitem[ReS2]{RS2}
\leavevmode\vrule height 2pt depth -1.6pt width 23pt, {\em Methods of modern
  mathematical physics II. Fourier analysis, self-adjoint\-ness}.
\newblock Academic Press, New York etc., 1975.

\bibitem[RSe3]{RS4}
\leavevmode\vrule height 2pt depth -1.6pt width 23pt, {\em Methods of modern mathematical physics
  IV. Analysis of operators}.
\newblock Academic Press, New York etc., 1978.

\bibitem[RoS]{RS}
{\sc Rothschild, L.P. {\rm and} Stein, E.M.}, Hypoelliptic differential
  operators and nilpotent groups.
\newblock {\em Acta Math.} {\bf 137} (1976),  247--320.

\bibitem[San]{SC}
{\sc S{\'a}nchez-Calle, A.}, Fundamental solutions and geometry of the sum of
  squares of vector fields.
\newblock {\em Invent.\ Math.} {\bf 78} (1984),  143--160.

\bibitem[Sik1]{Sik}
{\sc Sikora, A.}, Sharp pointwise estimates on heat kernels.
\newblock {\em Quart.\ J. Math.\ Oxford} {\bf 47} (1996),  371--382.

\bibitem[Sik2]{Sik3}
\leavevmode\vrule height 2pt depth -1.6pt width 23pt, Riesz transform, Gaussian
  bounds and the method of wave equation.
\newblock {\em Math.\ Z.} {\bf 247} (2004),  643--662.

\bibitem[Sil]{Sil}
{\sc Silverstein, M.L.}, {\em Symmetric Markov processes}.
\newblock Lect.\ Notes in Math. 426. Springer Verlag, Berlin, 1974.

\bibitem[Sim1]{bSim4}
{\sc Simon, B.}, Lower semicontinuity of positive quadratic forms.
\newblock {\em Proc.\ Roy.\ Soc.\ Edinburgh Sect.\ A} {\bf 79} (1977),
  267--273.

\bibitem[Sim2]{bSim5}
\leavevmode\vrule height 2pt depth -1.6pt width 23pt, A canonical decomposition
  for quadratic forms with applications to monotone convergence theorems.
\newblock {\em J.\ Funct.\ Anal.} {\bf 28} (1978),  377--385.

\bibitem[StW]{SW}
{\sc Stein, E.M. {\rm and} Weiss, G.}, {\em Introduction to Fourier analysis
  on Euclidian spaces}.
\newblock Princeton Mathematical Series 32. Princeton University Press,
  Princeton, New Jersey, 1971.

\bibitem[Stri]{Stri}
{\sc Strichartz, R.S.}, Multipliers on fractional Sobolev spaces.
\newblock {\em J. Math.\ Mech.} {\bf 16} (1967),  1031--1060.

\bibitem[Stro1]{Stroock3}
{\sc Stroock, D.W.}, Some applications of stochastic calculus to partial
  differential equations.
\newblock In {\em Eleventh Saint Flour probability summer school---1981 (Saint
  Flour, 1981)}, vol.\ 976 of Lecture Notes in Math.,  267--382. Springer,
  Berlin, 1983.

\bibitem[Stro2]{Stroock1}
\leavevmode\vrule height 2pt depth -1.6pt width 23pt, Diffusion semigroups corresponding to uniformly elliptic
  divergence form operators.
\newblock In {\sc Az{\'e}ma, J., Meyer, P.A. {\rm and} Yor, M.}, eds., {\em
  S{\'e}minaire de probabilit{\'e}s XXII}, Lecture Notes in Mathematics 1321.
  Springer-Verlag, Berlin etc., 1988,  316--347.

\bibitem[Stro3]{Stroock2}
\leavevmode\vrule height 2pt depth -1.6pt width 23pt, Estimates for the heat kernel of second order elliptic
  operators.
\newblock In {\em Nonlinear partial differential equations and their
  applications. Coll\`ege de France Seminar, Vol.\ XII (Paris, 1991--1993)},
  vol.\ 302 of Pitman Res.\ Notes Math.\ Ser.,  226--235. Longman Sci.\ Tech.,
  Harlow, 1994.

\bibitem[StV]{SV1}
{\sc Stroock, D.W. {\rm and} Varadhan, S. R.S.}, {\em Multidimensional
  diffusion processes}, vol.\ 233 of Grundlehren der Mathematischen
  Wissenschaften.
\newblock Springer-Verlag, Berlin, 1979.

\bibitem[Stu1]{Stu3}
{\sc Sturm, K.-T.}, Analysis on local Dirichlet spaces. I. Recurrence,
  conservativeness and $L_p$--Liouville properties.
\newblock {\em J. Reine Angew.\ Math.} {\bf 456} (1994),  173--196.

\bibitem[Stu2]{Stu4}
\leavevmode\vrule height 2pt depth -1.6pt width 23pt, Analysis on local
  Dirichlet spaces. II. Upper Gaussian estimates for the fundamental solutions
  of parabolic equations.
\newblock {\em Osaka J. Math.} {\bf 32} (1995),  275--312.

\bibitem[Stu3]{Stu5}
\leavevmode\vrule height 2pt depth -1.6pt width 23pt, Analysis on local
  Dirichlet spaces. III. The parabolic Harnack inequality.
\newblock {\em J. Math.\ Pures Appl.} {\bf 75} (1996),  273--297.

\bibitem[Stu4]{Stu2}
\leavevmode\vrule height 2pt depth -1.6pt width 23pt, The geometric aspect of
  Dirichlet forms.
\newblock In {\em New directions in Dirichlet forms}, vol.\ 8 of AMS/IP Stud.\
  Adv.\ Math.,  233--277. Amer.\ Math.\ Soc., Providence, RI, 1998.

\bibitem[Tay]{Tay3}
{\sc Taylor, M.E.}, {\em Partial differential equations. I{}. Basic theory},
  vol.\ 115 of Applied Mathematical Sciences.
\newblock Springer-Verlag, New York, 1996.

\bibitem[Tre]{Tre1}
{\sc Treves, F.}, {\em Introduction to pseudodifferential and Fourier integral
  operators}, vol.\ 1.
\newblock Plenum Press, New York etc., 1980.

\bibitem[Tru]{Tru2}
{\sc Trudinger, N.S.}, Linear elliptic operators with measurable coefficients.
\newblock {\em Ann.\ Scuola Norm.\ Sup.\ Pisa} {\bf 27} (1973),  265--308.

\bibitem[Var]{Varo10}
{\sc Varopoulos, N.T.}, Semigroup commutators under differences.
\newblock {\em Rev.\ Mat.\ Iberoamericana} {\bf 8} (1992),  1--43.

\end{thebibliography}
\end{document}